\documentclass[11pt]{article}
\usepackage{graphicx}
\usepackage{color}
\usepackage{amsmath}
\usepackage{amssymb}
\usepackage{amscd}
\usepackage{bbm}
\newcommand{\R}{\mathbb{R}}

\newcommand{\E}{\mathbb{E}}

\newcommand{\eps}{\varepsilon}

\newtheorem{Theorem}{Theorem}[section]
\newtheorem{Lemma}[Theorem]{Lemma}

\newtheorem{Definition}[Theorem]{Definition}

\newtheorem{Corollary}[Theorem]{Corollary}
\newtheorem{Remark}[Theorem]{Remark}

\numberwithin{equation}{section}

\def \proof {\noindent {\bf Proof.}\ \ }

\def \endproof
{{\mbox{}\nolinebreak\hfill\rule{2mm}{2mm}\par\medbreak}}
\def\IND{\mathbbm{1}}

\begin{document}
\title{On aggregation for heavy-tailed classes}
\author{Shahar Mendelson\footnote{Department of
Mathematics, Technion, I.I.T, Haifa 32000, Israel
\newline
{\sf email: shahar@tx.technion.ac.il}
\newline
Supported in part by the Mathematical Sciences Institute,
The Australian National University, Canberra, ACT 2601,
Australia.
 } }

\medskip
\maketitle
\begin{abstract}
We introduce an alternative to the notion of `fast rate' in Learning Theory, which coincides with the optimal error rate when the given class happens to be convex and regular in some sense.
While it is well known that such a rate cannot always be attained by a learning procedure (i.e., a procedure that selects a function in the given class), we introduce an aggregation procedure that attains that rate under rather minimal assumptions -- for example, that the $L_q$ and $L_2$ norms are equivalent on the linear span of the class for some $q>2$, and the target random variable is square-integrable.
\end{abstract}

\maketitle

\section{Introduction}
The focus of this article is on the question of {\it Prediction}: let ${\cal F}$ be a class of functions defined on a probability space $(\Omega,\mu)$ and let $X$ be distributed according to $\mu$. Given an unknown target random variable $Y$, one would like to find some $f \in {\cal F}$ for which, on average, predicting $f(X)$ instead of $Y$ is the most `cost effective'. If the pointwise cost is measured according to the squared loss, that is, if the price of predicting $f(X)$ instead of $Y$ is $(f(X)-Y)^2$, the goal is to identify, or at least approximate in some sense, the behaviour of the function that minimizes in ${\cal F}$ the {\it risk} $\E (f(X)-Y)^2$, where the expectation is taken with respect to the joint distribution of $X$ and $Y$ on $\Omega \times \R$. With that in mind, set
$$
f^*={\rm argmin}_{f \in {\cal F}} \E(f(X)-Y)^2,
$$
and assume, for the sake of simplicity, that the minimizer exists.

It should also be noted that there are other reasonable choices for the pointwise cost of predicting $f(X)$ instead of $Y$, and although our results are presented only for the squared loss, they may be extended to other convex loss functions, following the path of \cite{Men-LWCG}.

Unlike standard questions in Approximation Theory, in the prediction framework one has limited information: a random sample $(X_i,Y_i)_{i=1}^N$, selected independently according to the joint distribution of $X$ and $Y$. The hope is that a typical sample may be used to produce a (random) function in ${\cal F}$ that has almost the same `predictive capabilities' as the minimizer $f^*$.

\begin{Definition} \label{def:prediction}
For every integer $N$ and a base class ${\cal F}$, a {\it learning procedure}\footnote{This is sometimes called a {\it proper} learning procedure.} is a function $\Psi:(\Omega \times \R)^N \to {\cal F}$.

Setting $\tilde{f}=\Psi((X_i,Y_i)_{i=1}^N)$, and given $0<\delta<1$ and a set of potential targets ${\cal Y}$, the learning procedure $\Psi$ performs with an {\it error rate} of ${\cal E}_p({\cal F},N,\delta)$ if for every reasonable class of functions ${\cal F} \subset L_2(\mu)$ and $Y \in {\cal Y}$,
$$
\E \left((\tilde{f}(X)-Y)^2|(X_i,Y_i)_{i=1}^N\right) \leq \E(f^*(X)-Y)^2 + {\cal E}_p(F,N,\delta)
$$
with probability at least $1-\delta$ relative to the samples $(X_i,Y_i)_{i=1}^N$ (i.e., with respect to the $N$-product of the joint distribution of $X$ and $Y$ endowed on $(\Omega \times \R)^N$).
\end{Definition}
One would like to identify the `best' learning procedure $\Psi$, in the sense that the error rate ${\cal E}_p$ is as small as possible, find which features of ${\cal F}$ and ${\cal Y}$ govern ${\cal E}_p$, and study the way in which ${\cal E}_p$ scales with the sample size $N$.

\vskip0.3cm

Although it is not obvious from Definition \ref{def:prediction}, the effect the set of admissible targets ${\cal Y}$ has on the error rate ${\cal E}_p$ is rather small. In standard scenarios, ${\cal Y}$  consists of all random variables that are bounded by $1$ or, alternatively, that have rapidly decaying tails (e.g. - subgaussian or subexponential). However, as will be explained later, this type of condition can be relaxed considerably, and ${\cal Y}$ may be as large as the $L_2$ unit ball on the underlying probability space, rather than the $L_\infty$ one.

\subsection{Fast and slow rates} \label{sec:fast-slow-rates}
One frequently encounters in literature the terms `fast rate' and `slow rate', used to describe the behaviour of a learning procedure as a function of the sample size $N$. Unfortunately, the meaning of the two is somewhat ambiguous, and is often misinterpreted.

A common misapprehension is that `fast rate' means that ${\cal E}_p$ scales as $1/N$, and that a `slow rate' implies that ${\cal E}_p$ is of the order of $1/\sqrt{N}$; in reality, the situation is different. Indeed, on one hand, it is straightforward to construct examples of classes that are simply too rich for a rate of $1/N$ (or even of $1/\sqrt{N}$, for that matter), even in the realizable case, when $Y \in {\cal F}$; on the other, the `size' of ${\cal F}$ does not capture the correct behaviour of ${\cal E}_p$: if ${\cal F}=\{f_1,f_2\}$ and $Y$ happens to be a $1/\sqrt{N}$ perturbation of the mid-point $(f_1+f_2)/2$, no learning procedure can achieve an error rate that is better than $c/\sqrt{N}$ with probability at least $3/4$ using $N$ sample points and for a suitable absolute constant $c$ (see, e.g., \cite{AB} for a more precise statement).

Thus, a reasonable definition of the terms `fast rate' and `slow rate' must reflect the fact that the error rate is highly affected by the `location' of the target, as well as by the `complexity' of ${\cal F}$.

To avoid potential ambiguity, we will refrain from using the terms `fast rate' and `slow rate' in what follows. Instead, we will adopt the notion of {\it `optimistic rate'}, which is, roughly put, the rate one encounters when the location of the target is favourable, and should be considered as a more accurate version of the intuitive `fast rate' (see Section \ref{sec:opt} for the definition). For example, if ${\cal F}=\{f_1,f_2\}$, the optimistic rate is of the order of $1/N$ rather than $1/\sqrt{N}$, seemingly ignoring the possibility that $Y$ is a perturbation of the mid-point $(f_1+f_2)/2$ as above.

Note that this example shows that the optimistic rate may be, at times, unreachable by {\it any} learning procedure. Hence, if there is any hope of constructing a procedure that always attains the optimistic rate regardless of the location of the target, that procedure must be allowed the flexibility of selecting functions that are outside the given base class ${\cal F}$. Such procedures belong to the {\it model selection aggregation framework}.
\begin{Definition} \label{def:agg}
For an integer $N$, an aggregation procedure is a map $\Psi:(\Omega \times \R)^N \to L_2(\mu)$. The procedure has an error rate of ${\cal E}_p^{\rm agg}({\cal F},N,\delta)$ if for every reasonable class of functions ${\cal F}$ and every target $Y \in {\cal Y}$,
$$
\E \left((\tilde{f}(X)-Y)^2|(X_i,Y_i)_{i=1}^N\right) \leq \E(f^*(X)-Y)^2 + {\cal E}_p^{\rm agg}({\cal F},N,\delta)
$$
with probability at least $1-\delta$ relative to the $N$-product of the joint distribution of $X$ and $Y$, and for $\tilde{f}=\Psi((X_i,Y_i)_{i=1}^N)$.
\end{Definition}
Detailed surveys on the aggregation framework in a broad context may be found in \cite{Tsy,Cat,Nem}.

\vskip0.4cm
Thus, rather than restricting one to a learning procedure, i.e., forcing one to select functions from ${\cal F}$, the goal here is to construct an aggregation procedure that attains the optimistic rate under minimal assumptions on ${\cal F}$ and ${\cal Y}$.
\vskip0.4cm
The analysis of the aggregation procedure we will introduce below requires the use of some auxiliary classes that are connected to the given base class $F$; those will be denoted by $U,V$ and $H$. To avoid confusion, in what follows we will denote `generic' function classes by ${\cal F}$ and ${\cal K}$.

\subsection{The optimistic rate} \label{sec:opt}
The definition of the optimistic rate is based on the method developed in \cite{Men-LWC,Men-LWCG} for the analysis of the {\it Empirical Risk Minimization} procedure (ERM). We will outline the essentials of this method in what follows, but refer the reader to \cite{Men-LWC,Men-LWCG} for a more detailed description of the parameters involved, their role in the analysis of ERM and the way in which they may be computed in specific applications.

\begin{Definition} \label{def:ERM}
Given a sample $(X_i,Y_i)_{i=1}^N$ and a base class ${\cal F}$, the empirical minimizer in ${\cal F}$ is
$$
\hat{f} \in {\rm argmin}_{f \in {\cal F}} \frac{1}{N}\sum_{i=1}^N (f(X_i)-Y_i)^2,
$$
assuming, of course, that a minimizer exists.
\end{Definition}

From here on we will denote by $P_N h$ the empirical mean $\frac{1}{N}\sum_{i=1}^N h(X_i,Y_i)$.

\vskip0.3cm

Recall that $f^*={\rm argmin}_{f \in {\cal F}} \E(f(X)-Y)^2$, consider the squared excess loss functional relative to ${\cal F}$ and $Y$,
$$
{\cal L}_f(X,Y) = (f(X)-Y)^2-(f^*(X)-Y)^2,
$$
and observe that the minimizer in ${\cal F}$ of $P_N (f-Y)^2$ is also a minimizer in ${\cal F}$ of $P_N {\cal L}_f$. Thus, $P_N {\cal L}_{\hat f} \leq 0$, simply because ${\cal L}_{f^*}=0$, and, in particular,
$$
\hat{f} \in \{f \in {\cal F} : P_N {\cal L}_f \leq 0\}.
$$
It follows that if $(X_i,Y_i)_{i=1}^N$ is a sample for which
\begin{equation} \label{eq:inclusion}
\{f \in {\cal F}: \E {\cal L}_f \geq \eta\} \subset \{f \in {\cal F} : P_N {\cal L}_f >0\},
\end{equation}
then
$$
\E\left(( \hat{f}(X)-Y)^2 | (X_i,Y_i)_{i=1}^N\right) \leq \E (f^*(X)-Y)^2 + \eta,
$$
which is the type of result one is looking for.
\vskip0.3cm
To obtain \eqref{eq:inclusion}, note that for every $f \in {\cal F}$ and every sample $(X_i,Y_i)_{i=1}^N$,
\begin{align} \label{eq:excess-loss-expansion}
P_N {\cal L}_f & \geq \frac{1}{N} \sum_{i=1}^N (f-f^*)^2(X_i) + 2\E (f^*(X)-Y)(f-f^*)(X)
\\
& -2 \left|\frac{1}{N}\sum_{i=1}^N (f^*(X_i)-Y_i)(f-f^*)(X_i) - \E (f^*(X)-Y)(f-f^*)(X)\right|, \nonumber
\end{align}
i.e., $P_N {\cal L}_f$ is lower bounded by a sum of (random) quadratic and multiplier components, and a deterministic term,  $2\E (f^*(X)-Y)(f-f^*)(X)$, which calibrates the `location' of the target $Y$ relative to ${\cal F}$.

The optimistic rate is defined based on the belief that the location of $Y$ is favourable in the sense that for every $f \in {\cal F}$,
\begin{equation} \label{eq:B-condition}
\E (f^*(X)-Y)(f-f^*)(X) \geq 0.
\end{equation}
It is straightforward to verify that \eqref{eq:B-condition} is satisfied in two important cases. Firstly, when ${\cal F} \subset L_2$ happens to be closed and convex, in which case, \eqref{eq:B-condition} follows from the characterization of the metric projection onto a closed, convex set in an inner-product space. Secondly, for an arbitrary class ${\cal F}$ and a target $Y=f^*(X)+\xi$, where $f^* \in {\cal F}$ and $\xi$ is mean-zero and independent of $X$.

\vskip0.3cm

\begin{Definition} \label{def:small-ball}
A class ${\cal F}$ satisfies a small-ball condition with constants $\kappa_0$ and $\eps$, if for every $f_1,f_2 \in {\cal F} \cup\{0\}$,
\begin{equation} \label{eq:small-ball}
Pr(|f_1-f_2| \geq \kappa_0\|f_1-f_2\|_{L_2}) \geq \eps.
\end{equation}
\end{Definition}
The small-ball condition is a rather minimal assumption on ${\cal F}$ -- it is a uniform lower estimate on the probability that $|f_1-f_2|/\|f_1-f_2\|_{L_2}$ is sufficiently far from zero for every pair of distinct functions $f_1,f_2 \in {\cal F} \cup \{0\}$.

One may find in \cite{Men-LWC,Men-LWCG} several examples of classes that satisfy a small-ball condition. For our purposes, the most significant example is when $q>2$ and the $L_q$ and $L_2$ norms are $L$-equivalent on ${\cal F}$, in the sense that for every $f_1,f_2 \in {\cal F} \cup \{0\}$, $\|f_1-f_2\|_{L_q} \leq L \|f_1-f_2\|_{L_2}$. In such a case, the Paley-Zygmund inequality \cite{dlPG} shows that \eqref{eq:small-ball} holds for constants $\kappa_0$ and $\eps$ that depend only on $q$ and $L$.

\vskip0.4cm

Let $D$ be the unit ball in $L_2(\mu)$. Given a class of functions ${\cal F} \subset L_2(\mu)$, set $\{G_f : f \in {\cal F}\}$ to be the canonical gaussian process indexed by ${\cal F}$ and put
$$
\E\|G\|_{\cal F} = \sup\{ \E \sup_{ f \in {\cal F}^\prime} G_f : {\cal F}^\prime \subset {\cal F}, \ {\cal F}^\prime {\rm \ is \ finite \ } \}.
$$

\begin{Definition} \label{def:fixed-point}
 For ${\cal F} \subset L_2$, let ${\rm star}({\cal F})=\{\lambda f : \ 0 \leq \lambda \leq 1, \ f \in {\cal F}\}$ be
the star-shaped hull of ${\cal F}$ around $0$, and set ${\cal F}-{\cal F} =\{f-h : f,h \in {\cal F}\}$.
Let
$$
U=\left\{\frac{f_1+f_2}{2} : f_1,f_2 \in {\cal F} \right\}
$$
and set $H={\rm star}(U-U)$.

\vskip0.3cm
Finally, for $\zeta>0$, let
$$
r_{Q,1}({\cal F},\zeta)=\inf\left\{r>0 \ : \ \E\|G\|_{(H-H) \cap r D} \leq \zeta r \sqrt{N} \right\},
$$
and
$$
r_{Q,2}({\cal F},\zeta)=\inf\left\{r>0 \ : \E \sup_{w \in (H-H) \cap r D} \left|\frac{1}{\sqrt{N}} \sum_{i=1}^N \eps_i w(X_i) \right| \leq \zeta r \sqrt{N} \right\},
$$
where $(\eps_i)_{i=1}^N$ are independent, symmetric $\{-1,1\}$-valued random variables that are independent of $(X_i)_{i=1}^N$, and the expectation is taken with respect to both $(\eps_i)_{i=1}^N$ and $(X_i)_{i=1}^N$.
\end{Definition}

Note that ${U}$ is only slightly richer than ${\cal F}$: it contains ${\cal F}$ and all the midpoints of intervals whose ends belong to ${\cal F}$. If ${\cal F}$ happens to be convex, then ${U}={\cal F}$, but in general, ${U}$ is much smaller than the convex hull of ${\cal F}$. Also, $H={\rm star}(U-U)$ is star-shaped around $0$, centrally symmetric, and contains ${\cal F}-{\cal F}$; hence, both ${\cal F}$ and ${\cal F}-{\cal F}$ belong to $H-H$.

\vskip0.3cm

The parameters $r_{Q,1}$ and $r_{Q,2}$ measure the `local' complexity of the indexing class: from a statistical point of view, the two capture the correlation of the indexing class with various forms of random noise. The reader may find a more detailed explanation of their role in \cite{Men-LWC} and  \cite{Men-LWCG}.

It should be noted that the definitions of $r_{Q,1}$ and $r_{Q,2}$ in \cite{Men-LWCG} appear to be slightly different from the ones defined above. However, the reason for the difference is that in \cite{Men-LWCG} one considers a convex base class, while here ${\cal F}$ need not be convex. If ${\cal F}$ happens to be convex then $U={\cal F}$, $U-U={\cal F}-{\cal F}$ is convex and centrally symmetric, and $H={\cal F}-{\cal F}$; therefore $H-H=2H=2({\cal F}-{\cal F})$ and the definitions above coincide with the ones from \cite{Men-LWCG} up to a factor of $2$, which is only an issue of normalization.
\vskip0.4cm

The third and final complexity parameter is also a minor modification of a similar parameter from \cite{Men-LWC,Men-LWCG}. It will be used to study the multiplier component in the decomposition \eqref{eq:excess-loss-expansion} of the excess squared-loss functional.
\begin{Definition} \label{def:complexity-linear}
Let ${\cal F} \subset L_2$ be the given base class and set ${U}$ and $H$ as above. For every $u_0 \in {U}$ consider the random function
$$
\phi_{{\cal F},N,u_0}(r)=\frac{1}{\sqrt{N}} \sup_{\{w \in {\rm star}({U}-u_0) \cap rD\}} \left|\sum_{i=1}^N \eps_i (u_0(X_i)-Y_i)w(X_i)\right|
$$
and set
$$
r_{M}({\cal F},\zeta,\delta,u_0)=\inf\left\{ r > 0 : Pr \left(\phi_{{\cal F},N,u_0}(r) \leq  r^2 \zeta \sqrt{N} \right) \geq 1-\delta \right\}.
$$
\end{Definition}

The importance of $r_M$ and the way it may be used to upper bound the multiplier component can be seen in the next lemma from \cite{Men-LWCG}:
\begin{Lemma} \label{lemma:multi-and-zero-level-sym}
Let $u_0 \in {U}$, $0<\delta<1$, $\kappa>0$ and set $r =2 r_M({\cal F},\kappa/4,\delta/2,u_0)$. Put $\xi=u_0(X)-Y$, and given a sample $(X_i,Y_i)_{i=1}^N$ set $\xi_i=u_0(X_i)-Y_i$. Then, with probability at least $1-\delta$, for every $u \in {U}$ that satisfies $\|u-u_0\|_{L_2} \geq r$, one has
$$
\left|\frac{1}{N}\sum_{i=1}^N \xi_i(u-u_0)(X_i) - \E \xi(u-u_0)(X) \right| \leq \kappa \max\{\|u-u_0\|_{L_2}^2,r^2\}.
$$
\end{Lemma}

With all the complexity terms in place, one may derive an error estimate for ERM, performed in any subset of $U$ that satisfies the `optimistic' assumption. Indeed, the following is a minor modification of Lemma 5.2 from \cite{Men-LWCG}, originally formulated for a convex base class, though it is straightforward to verify that the convexity condition may be relaxed.

\vskip0.3cm

In the setup we are interested in, ${\cal F} \subset L_2$ is the given base class, $U$ and $H$ are defined as above and $H$ satisfies a small-ball condition with constants $\kappa_0$ and $\eps$. Fix $Y \in L_2$ and $V \subset U$, put $v^*={\rm argmin}_{v \in V} \|v-Y\|_{L_2}$ and let $\hat{v}$ be the empirical minimizer in $V$ of the squared loss functional.

\begin{Theorem} \label{thm:error-rate-convex}
There exists an absolute constant $c_0$ and for every $\kappa_0>0$ and $0<\eps<1$ there exist constants $c_1,c_2$ and $c_3$ that depend only on $\kappa_0$ and $\eps$ for which the following holds.
 If for every $v \in V$, $\E (v^*(X)-Y)(v-v^*)(X) \geq 0$, then probability at least $1-\delta-2\exp(-c_0\eps^2 N)$,
$$
\E \left((\hat{v}(X)-Y)^2 | (X_i,Y_i)_{i=1}^N\right) \leq \E(v^*(X)-Y)^2) + r^2(v^*),
$$
where
\begin{equation*}
r(v^*) = 2\max\left\{r_M({\cal F},c_1,\delta/2,v^*),r_{Q,1}({\cal F},c_2),r_{Q,2}({\cal F},c_3)\right\}.
\end{equation*}
\end{Theorem}

Since $v^*$ is not known, one has to use a uniform version of $r(v^*)$ as a complexity parameter. This uniform version is the optimistic rate:
\begin{Definition} \label{def:optimistic-rate}
Given a base class $F$, the optimistic rate in $F$ is defined by
\begin{equation} \label{eq:r-opt}
r_{{\rm opt}}(F,\delta,c_1,c_2,c_3) = 2 \sup_{u_0 \in {U}} \max\left\{r_M^2(F,c_1,\delta/2,u_0),r_{Q,1}^2({F},c_2),r_{Q,2}^2({F},c_3)\right\}.
\end{equation}
\end{Definition}
In what follows, $H={\rm star}(U-U)$ will satisfy the small-ball condition with constants $\kappa_0$ and $\eps$, and $c_1,c_2$ and $c_3$ will be chosen as constants that depend only on $\kappa_0$ and $\eps$. To avoid cumbersome notation, we will not specify in what follows that $r_{{\rm opt}}$ depends on $F$, $\delta$ and the constants $c_1,c_2,c_3$, but their choice will be made clear.

\vskip0.4cm
When $F$ happens to be convex, $U=F$, and upon selecting $V=F$ it follows that $\E (f^*(X)-Y)(f-f^*)(X) \geq 0$ for every $f \in F$. Thus, Theorem \ref{thm:error-rate-convex} extends the main result from \cite{Men-LWC} on the performance of ERM in a convex class and relative to the squared loss. Moreover, since the `optimistic assumption' includes the choice of $Y=f^*(X)+\xi$ for $f^* \in F$ and $\xi$ that is independent of $X$, the results from \cite{LM-subgaussian,Men-LWC} indicate that $r_{{\rm opt}}^2$ captures the {\it minimax rate}\footnote{Roughly put, the minimax rate is the best possible error rate one may achieve by {\it any} learning procedure, i.e., by any $\Psi:(\Omega \times \R)^N \to F$.} in $F$ under mild structural assumptions on that class.
\vskip0.4cm
Therefore, the optimistic rate $r_{\rm opt}$ is defined as what is essentially the best possible rate that any learning procedure may achieve in $F$ when the target $Y$ is in a `good location' relative to $F$ in the sense of \eqref{eq:B-condition}. And, when $F$ happens to be convex, every target $Y \in L_2$ is in a `good location'.

Having said that, let us emphasize once again that the problem we wish to address occurs when the location of the target is less favourable, and in which case no learning procedure can achieve the optimistic rate.

\vskip0.3cm

We will show that there is an {\it aggregation procedure} that always achieves the optimistic rate when the $L_2$ and $L_q$ norms are equivalent on ${\rm span}(F)$ for some $q>2$, and $Y \in L_2$ -- thus overcoming the possible problem that may occur when the target $Y$ is not in a `good location' relative to the given class. Let us stress that what allows one to attain the optimistic rate is that the procedure used in an aggregation procedure, and thus may take values in $L_2(\mu)$, rather than a learning procedure, which is restricted to values in $F$.

\vskip0.3cm
\begin{Theorem} \label{thm:main}
For every $L \geq 1$ and $q>2$ there are constants $c_0,c_1,c_2$ and $c_3$ that depend only on $L$ and $q$ for which the following holds. Let $F \subset L_2$ be the given base class and let ${U}$ and $H$ be as above. Assume that for every $w \in H-H$, $\|w\|_{L_q} \leq L \|w\|_{L_2}$. Then, there is an aggregation procedure $\Psi:(\Omega \times \R)^N \to L_2(\mu)$ for which, for every $Y \in L_2$, with probability at least $1-\delta-2\exp(-c_0 N)$,
$$
\E (\tilde{f}(X)-Y)^2 | (X_i,Y_i)_{i=1}^N) \leq \E(f^*(X)-Y)^2) + r_{{\rm opt}},
$$
where $\tilde{f}=\Psi((X_i,Y_i)_{i=1}^N)$ and
\begin{equation} \label{eq:error-estimate-in-main}
r_{{\rm opt}}=  2\sup_{u_0 \in {U}} \max\left\{r_M^2(F,c_1,\delta/4,u_0),r_{Q,1}^2(F,c_2),r_{Q,2}^2(F,c_3)\right\}.
\end{equation}
\end{Theorem}

To put Theorem \ref{thm:main} is some perspective, note that in the standard framework of aggregation, $F$ is a finite dictionary and both the dictionary and the target are bounded in $L_\infty$ (see, for example, \cite{Tsy,Lec-HDR} and references therein). Within that framework one has the following:

\begin{Theorem} \label{thm:LM-agg} \cite{LM-agg}
There exists an aggregation procedure $\Psi$ for which the following holds.
Assume that $F$ is a finite dictionary consisting of functions that are bounded by $1$, and assume that the target $Y$ is bounded by $1$ as well. Then, for every $x>0$, with probability at least $1-2\exp(-x)$,
$$
\E\left((\tilde{f}(X)-Y)^2|(X_i,Y_i)_{i=1}^N\right) \leq \E(f^*(X)-Y)^2 + c(1+x)\frac{\log |F|}{N},
$$
where $\tilde{f}=\Psi((X_i,Y_i)_{i=1}^N)$.
\end{Theorem}

In comparison, when applied to a finite dictionary, Theorem \ref{thm:main} leads to the following:

\begin{Corollary} \label{cor:subgaussian}
Let $F$ be a finite dictionary and set $H$ as above. Assume that $H-H$ is $L$-subgaussian, in the sense that for every $w_1,w_2 \in H-H$, and every $p \geq 2$, $\|w_1-w_2\|_{L_p} \leq L\sqrt{p}\|w_1-w_2\|_{L_2}$. If $Y \in L_q$ for some $q>2$ then with probability at least $1-\delta$,
\begin{equation} \label{eq:example-error-rate}
\E (\tilde{f}(X)-Y)^2 | (X_i,Y_i)_{i=1}^N) \leq \E(f^*(X)-Y)^2)+c_1 \|f^*-Y\|_{L_q}^2 \frac{\log |F|}{N},
\end{equation}
where $c_1$ depend only on $q$, $L$ and $\delta$.
\end{Corollary}
The proof of Corollary \ref{cor:subgaussian} will be presented in Section \ref{sec:finite-dictionary}. It is well known that the best error rate a learning procedure may attain for a finite dictionary is of the order of $\sqrt{{\log |F|}/{N}}$ (see, e.g. \cite{Tsy,Lec-HDR}), which is not remotely close to $r_{\rm opt}$, as the latter is of the order of $(\log |F|)/N$. Moreover, the error rate in \eqref{eq:example-error-rate} scales well with  $\|f^*-Y\|_{L_q}$: it tends to zero when $Y$ approaches $F$ and the problem becomes `more realizable', in which case one expects a zero-error when $N \geq c\log |F|$.

\vskip0.4cm

The aggregation procedure we will introduce here is a `close family member' of the one from \cite{LM-agg}, but with many significant and unavoidable changes. It should be noted that Audibert obtained in \cite{Aud1} the same estimate as in Theorem \ref{thm:LM-agg}, but using a different aggregation procedure -- the {\it empirical star algorithm}, and it is not clear whether it is possible to obtain a version of Theorem \ref{thm:main} using an analog of the empirical star algorithm. Moreover, the empirical star algorithm involves running ERM on the star-hull of $F$ and the empirical minimizer; therefore, one has to apply ERM to an infinite class even if the dictionary is finite. In contrast, the procedure suggested here uses ERM on a well-chosen $V \subset U$; hence, if $F$ is finite, so is $V$.

Unlike the bounded case, aggregation in unbounded situations was not fully understood. The benchmark result in that direction is due to Audibert \cite{Aud2} and independently to Juditsky, Rigollet and Tsybakov \cite{JRT}, who obtained the following estimate on the expected risk when the class is bounded but the target may be unbounded:
\begin{Theorem} \label{thm:Aud-unbounded}
There is an aggregation procedure $\Psi:(\Omega \times \R)^N \to L_2(\mu)$ for which the following holds.
Assume that $F$ is a finite dictionary consisting of functions bounded by $1$ and assume that $Y \in L_q$ for $q \geq 2$. Then setting $\tilde{f}=\Psi((X_i,Y_i)_{i=1}^N)$,
$$
\E \left(\E\left((\tilde{f}(X)-Y)^2|(X_i,Y_i)_{i=1}^N\right)\right) \leq \E(f^*(X)-Y)^2 + C(q,\|Y\|_{L_q})\left(\frac{\log |F|}{N}\right)^{2/(q+2)};
$$
moreover, this estimate is optimal -- up to the constant $C$.
\end{Theorem}

It is interesting to see the subtle differences between the assumptions used in Theorem \ref{thm:Aud-unbounded} and the ones from Corollary \ref{cor:subgaussian}. In the former, the class is assumed to be bounded in $L_\infty$, while in the latter, the class is $L$-subgaussian, which is a different type of condition: it implies {\it norm equivalence} rather than having a bounded diameter with respect to some (possibly strong) norm.

As noted in \cite{LM-lin-agg}, statistical procedures may behave in a very different way when one assumes even a weak norm equivalence rather than an $L_\infty$ bound, and the same phenomenon is true here as well: although the dictionary may consist of unbounded functions, the norm equivalence gives sufficient information to ensure an error rate of $N^{-1}\log |F|$ rather than much slower $(N^{-1}\log |F|)^{2/(q+2)}$. Moreover, the error rates in Theorem \ref{thm:LM-agg} and Theorem \ref{thm:Aud-unbounded}  do not scale well with the distance between $Y$ and $F$ and do not improve even when the problem is arbitrarily close to being realizable.

\vskip0.4cm
We end this introduction with some notation. Throughout, absolute constants are denoted by $c,c_1...$, etc. Their value may change from line to line. When a constant depends on a parameter $\alpha$ it will be denoted by $c(\alpha)$. $A \lesssim B$ means that $A \leq cB$ for an absolute constant $c$, and $A \lesssim_\alpha B$ implies that the constant depends on the parameter $\alpha$. The analogous two-sided inequalities are denoted by $A \sim B$ and $A \sim_\alpha B$.

For a set $A$, let $\IND_A$ be its indicator function and put $|A|$ to be its cardinality.

Finally, let us mention that we will abuse notation and write $\|z\|_{L_2}$ for the $L_2$ norm of the function $z$, without specifying the exact probability space on which the integration is performed. For example, $\|f-Y\|_{L_2}^2=\E(f(X)-Y)^2$, while $\|f-f^*\|_{L_2}^2=\E (f-f^*)^2(X)$. We will denote the unit ball in $L_2$ by $D$ and the unit sphere by $S(L_2)$, again, without specifying the underlying space.

\section{The aggregation procedure} \label{sec:the-procedure}
The aggregation procedure presented here follows the general path of \cite{LM-agg} -- though with many essential modifications. The core difference between the method of proof used in \cite{LM-agg} and the one we employ here is unavoidable, as the former is based on a two-sided concentration estimate on empirical means which is simply false for heavy-tailed functions. Most notably,  two-sided empirical estimates on $L_2$ distances play a central role in \cite{LM-agg} and one has to find an alternative to these concentration-based bounds. To that end, we will introduce an  empirical `isomorphic' upper estimate on $L_2$ distances, which is based on the idea of median-of-means, and which will be complemented by an `almost isometric' lower bound. Both bounds hold for any two class-members that are not `very close', and under a weak moment assumption: that for some $q>2$ the $L_q$ and $L_2$ norms are $L$-equivalent on the class. These results are of independent interest and are likely to have many other applications.

The accurate formulation and proof of the `isomorphic' estimate may be found in Section \ref{sec:median}, while the `almost isometric' lower bound is presented in Section \ref{sec:isometric}.

\vskip0.5cm
The aggregation procedure consists of two stages. Given a base class $F$, one must first identify a subset $V \subset F$, which is selected in a data-dependent way, and which consists of well-behaved functions in a sense that will be clarified below. Then, in the second stage, one applies ERM to the set of midpoints of pairs of elements in $V$ (a set which contains $V$ as well), i.e., to
$$
W = \left\{\frac{v_1+v_2}{2} : v_1,v_2 \in V\right\} \subset U,
$$
using a second, independent sample.
\vskip0.3cm
We begin with the following observation:

\begin{Lemma} \label{lemma:V-class}
Let $C \geq 1$, $r>0$ and $0<\theta \leq 1/32$, and consider $V \subset F$ that satisfies the following:
\begin{description}
\item{$\bullet$} $f^* \in V$ (where, as always, $f^*={\rm argmin}_{f \in F} \|f-Y\|_{L_2}$).
\item{$\bullet$} For every $v \in V$, $\|v-Y\|_{L_2}^2 \leq \|f^*-Y\|_{L_2}^2 +\max\{Cr^2,\theta {\rm diam}^2(V,L_2)\}$.
\end{description}
Let $W = \left\{(v_1+v_2)/{2} : v_1,v_2 \in V\right\}$, put $w^*={\rm argmin}_{w \in W} \|w-Y\|_{L_2}$ and set $\hat{w}$ to be the empirical minimizer in $W$.

If $(X_i,Y_i)_{i=1}^N$ is a sample for which, for every $w \in W$,
\begin{align} \label{eq:erm-in-W}
& \|w-Y\|_{L_2}^2 - \|w^*-Y\|_{L_2}^2
\\
\leq & \frac{1}{N}\sum_{i=1}^N \left((w(X_i)-Y_i)^2 - (w^*(X_i)-Y_i)^2\right) +  \max\left\{Cr^2,\theta \|w-w^*\|_{L_2}^2 \right\}, \nonumber
\end{align}
then
$$
\E \left((\hat{w}-Y)^2|(X_i,Y_i)_{i=1}^N \right) \leq \E(f^*(X)-Y)^2 + 2Cr^2.
$$
\end{Lemma}
Lemma \ref{lemma:V-class} implies that if $V$ consists of functions whose excess risk (relative to $F$) is either small ($\leq Cr^2$), or, alternatively, at least smaller than a fixed proportion of the square of the diameter of $V$, and if one is given a sample for which the oracle type inequality \eqref{eq:erm-in-W} holds, then ERM performed in $W$ using that sample selects a function whose excess risk is at most $2Cr^2$.
\vskip0.4cm
Naturally, at this point Lemma \ref{lemma:V-class} is somewhat speculative, as it contains two substantial `if's'. For the lemma to be of any use, one has to construct $V$ using a random sample and without knowing the identity of $f^*$, and then to establish the oracle type inequality in $W$ using a second, independent sample.
\vskip0.3cm

\proof
Set $d_V={\rm diam}(V,L_2)$, note that ${\rm diam}(W,L_2)=d_V$ and that
\begin{equation} \label{eq:err-w-in-proof}
\frac{1}{N}\sum_{i=1}^N (\hat{w}(X_i)-Y_i)^2 - \frac{1}{N} \sum_{i=1}^N (w^*(X_i)-Y_i)^2 \leq 0.
\end{equation}
Consider two cases: firstly, if $d_V \leq \sqrt{C}r$ then by \eqref{eq:erm-in-W}
$$
\E\left((\hat{w}(X)-Y)^2|(X_i,Y_i)_{i=1}^N\right) \leq \E (w^*(X)-Y)^2  +Cr^2 \leq \E (f^*(X)-Y)^2  +Cr^2
$$
because $f^* \in V \subset W$.

Secondly, assume that $d_V \geq \sqrt{C}r$. Since $f^* \in V$, there is some $v \in V$ for which $\|v-f^*\|_{L_2} \geq d_V/2$. Set $w=(v+f^*)/2 \in W$ and observe that by the uniform convexity of the $L_2$ norm and the definition of $V$,
\begin{align*}
\|w^*-Y\|_{L_2}^2 \leq & \|w-Y\|_{L_2}^2 =  \frac{1}{2}\|v-Y\|_{L_2}^2 + \frac{1}{2}\|f^*-Y\|_{L_2} - \frac{1}{4}\|v-f^*\|_{L_2}^2
\\
\leq & \|f^*-Y\|_{L_2}^2 +\max\left\{Cr^2,\theta d_V^2\right\} - \frac{d_V^2}{16}
\\
\leq & \|f^*-Y\|_{L_2}^2 + Cr^2 -\left(\frac{1}{16}-\theta\right)d_V^2.
\end{align*}
Combining this with \eqref{eq:erm-in-W} applied to $\hat{w}$, and with \eqref{eq:err-w-in-proof}, and recalling that $\theta \leq 1/32$,
\begin{align*}
& \E\left((\hat{w}(X)-Y)^2|(X_i,Y_i)_{i=1}^N\right)
\\
\leq & \E(f^*(X)-Y)^2 + Cr^2 + \theta d_V^2 + \left(\E(w^*(X)-Y)^2-\E(f^*(X)-Y)^2\right)
\\
\leq & \E(f^*(X)-Y)^2 + 2Cr^2 -\left(\frac{1}{16}-2\theta\right)d_V^2
\\
\leq & \E(f^*(X)-Y)^2 + 2Cr^2.
\end{align*}
\endproof

Next, we shall identify sufficient conditions that allow one to construct the set $V$ as in Lemma \ref{lemma:V-class} in a data-dependent way.
\subsection{The construction of $V$}
Given a class $F$, recall that ${U}=\left\{(f_1+f_2)/2 : f_1,f_2 \in F \right\}$.

What will assume the role of the empirical mean $\frac{1}{N}\sum_{i=1}^N (f-h)^2(X_i)$ as a way of estimating $L_2$ distances, is the following median-of-means functional, which is more stable than the empirical mean when dealing with heavy-tailed functions:
\begin{Definition}
Let $1 \leq \ell \leq N$ and set $I_j =\{\ell j +1,...,\ell(j+1)\} \subset \{1,...,N\}$ for $0 \leq j \leq \lfloor N/\ell\rfloor \equiv M-1$. For $v \in \R^N$ let
${\rm Med}_\ell(v)$ to be the median of the vector of means $(\ell^{-1}\sum_{i \in I_j} v_i)_{j=0}^{M-1} \in \R^M$.
\end{Definition}
Thus, $I_0,...,I_{M-1}$ are disjoint subsets of $\{1,...,N\}$, each of cardinality $\ell$, and ${\rm Med}_\ell(v)$ is the median of the means taken over the `blocks' $I_j$.

\begin{Definition} \label{def:A-f}
Fix $r_{U}>0$, $u_0 \in {U}$, $1 \leq \ell \leq N$, $0<\alpha \leq 1 \leq \beta$, and set
$\rho=\left({\alpha}/{20\beta}\right)^2 \leq {1}/{400}$.

Let ${\cal A}_{u_0}$ be the set of $N$-samples $(X_i,Y_i)_{i=1}^N$ for which the following holds:
\begin{description}
\item{$\bullet$} for every $u \in {U}$
\begin{align*}
& \left|\frac{1}{N}\sum_{i=1}^N (u_0(X_i)-Y_i)(u-u_0)(X_i)- \E (u_0(X)-Y)(u-u_0)(X)\right|
\\
\leq & \rho \max\left\{r_{U}^2, \|u-u_0\|_{L_2}^2\right\};
\end{align*}
\item{$\bullet$} if $u_1,u_2 \in {U}$ and $\|u_1-u_2\|_{L_2} \geq r_{{U}}$, then
$$
\left(\frac{1}{N}\sum_{i=1}^N (u_1-u_2)^2(X_i)\right) \geq (1-\rho)\|u_1-u_2\|_{L_2}^2;
$$
\item{$\bullet$} if $u_1,u_2 \in {U}$ and $\|u_1-u_2\|_{L_2} \geq r_{{U}}$, then
$$
\alpha\|u_1-u_2\|_{L_2} \leq {\rm Med}_\ell(|u_1-u_2|(X_i))_{i=1}^N \leq \beta \|u_1-u_2\|_{L_2},
$$
and if $\|u_1-u_2\|_{L_2} < r_{{U}}$ then
$$
{\rm Med}_\ell(|u_1-u_2|(X_i))_{i=1}^N \leq \beta r_{{U}}.
$$
\end{description}
\end{Definition}

At this point, the right choice of $r_{{U}}$, $\ell$, $\alpha$ and $\beta$ is not clear, nor that ${\cal A}_{u_0}$ is nonempty, for that matter.
\vskip0.4cm
At last, we are ready to define the aggregation procedure:
\begin{Definition} \label{def:sets-in-agg}
Let $(X_i,Y_i)_{i=1}^{2N}$ be a $2N$-sample and set ${\cal D}_1=(X_i,Y_i)_{i=1}^N$ and ${\cal D}_2=(X_i,Y_i)_{i=N+1}^{2N}$. Recall that
$$
\hat{f}={\rm argmin}_{f \in F} \frac{1}{N}\sum_{i=1}^N (f(X_i)-Y_i)^2,
$$
and let
\begin{align} \label{eq:def-hat-F}
V({\cal D}_1)=\Bigg\{f \in F : &  \frac{1}{N}\sum_{i=1}^N (f(X_i)-Y_i)^2 \leq  \frac{1}{N}\sum_{i=1}^N (\hat{f}(X_i)-Y_i)^2  \nonumber
\\
+ & 3\max\left\{ r_{{U}}^2, \rho\alpha^{-2}{\rm Med}_{\ell}^2\left(|\hat{f}-f|(X_i)\right)_{i=1}^N\right\}\Bigg\}.
\end{align}

Set
$$
W({\cal D}_1) = \left\{\frac{v_1+v_2}{2} : v_1,v_2 \in V({\cal D}_1) \right\}
$$
and define the aggregation procedure by
\begin{equation} \label{eq:procedure}
\tilde{w}={\rm argmin}_{w \in W({\cal D}_1)} \frac{1}{N}\sum_{i=N+1}^{2N} (w(X_i)-Y_i)^2.
\end{equation}
\end{Definition}

\begin{Theorem} \label{thm:procedure}
Let $w^*={\rm argmin}_{w \in V({\cal D}_1)} \E (w(X)-Y)^2$. If ${\cal D}_1 \in {\cal A}_{f^*}$ and ${\cal D}_2 \in {\cal A}_{w^*}$, then $V({\cal D}_1)$ and $W({\cal D}_1)$ satisfy the conditions of Lemma \ref{lemma:V-class} for $\theta<1/32$, $C=6$ and $r=r_{U}$. In particular, for such a $2N$-sample, and if $\tilde{w}$ is the empirical minimizer selected in $W({\cal D}_1)$ using the sample ${\cal D}_2$, one has
$$
\E\left((\tilde{w}(X)-Y)^2|{\cal D}_1\right) \leq \E\left(f^*(X)-Y\right)^2 + 6r_{{U}}^2.
$$
\end{Theorem}
Theorem \ref{thm:main} follows from Theorem \ref{thm:procedure}, once one shows that $r_U^2$ can be selected to be $r_{{\rm opt}}$ for the right choice of constants, and that the probability of the events ${\cal A}_{u_0}$ is sufficiently high for every $u_0 \in U$.
\vskip0.3cm

The rest of this section is devoted to the proof of Theorem \ref{thm:procedure}. The proof that each ${\cal A}_{u_0}$ is a large event will be presented in Section \ref{sec:the-event}.

\begin{Remark}
Note that if $F$ is finite, so is $W({\cal D}_1)$ -- which may be much smaller than $F$. Thus, if the original dictionary is finite, then unlike the empirical star algorithm, the second step in the aggregation procedure is carried out on a finite set.
\end{Remark}

\subsection{Proof of Theorem \ref{thm:procedure}}
Given a set ${\cal K} \subset {U}$, put $h^*={\rm argmin}_{h \in {\cal K}} \E (h(X)-Y)^2$ and set
$$
{\cal L}_h^{{\cal K}}(X,Y)=(h(X)-Y)^2-(h^*(X)-Y)^2
$$
to be the square excess loss functional associated with ${\cal K}$.
\begin{Lemma} \label{lemma:key-observation}
If $(X_i,Y_i)_{i=1}^N \in {\cal A}_{h^*}$ then for every $h \in {\cal K}$,
\begin{equation} \label{eq:est-implies-pre}
\E {\cal L}^{{\cal K}}_h \leq P_N {\cal L}^{{\cal K}}_h  + 3\max\{r_{{U}}^2, \rho\|h-h^*\|_{L_2}^2\}.
\end{equation}
\end{Lemma}

\proof
Set $\xi(X,Y)=h^*(X)-Y$, let $(X_i,Y_i)_{i=1}^N \in {\cal A}_{h^*}$ and put $\xi_i=h^*(X_i)-Y_i$. Note that for every $h \in {\cal K}$,
\begin{align*}
P_N {\cal L}^{{\cal K}}_h = & \frac{2}{N}\sum_{i=1}^N \xi_i (h-h^*)(X_i) + \frac{1}{N} \sum_{i=1}^N (h-h^*)^2(X_i)
\\
\geq & \E {\cal L}^{{\cal K}}_h - 2\left|\frac{1}{N} \sum_{i=1}^N \xi_i (h-h^*)(X_i) - \E \xi(h-h^*)(X)\right|
\\
+ & \frac{1}{N} \sum_{i=1}^N (h-h^*)^2(X_i) - \|h-h^*\|_{L_2}^2.
\end{align*}
Recall that if $(X_i,Y_i)_{i=1}^N \in {\cal A}_{h^*}$ and $\|h-h^*\|_{L_2} \geq r_{{U}}$, one has
$$
\frac{1}{N} \sum_{i=1}^N (h-h^*)^2(X_i) \geq (1-\rho) \|h-h^*\|_{L_2}^2,
$$
and
$$
\left|\frac{1}{N} \sum_{i=1}^N \xi_i (h-h^*)(X_i) - \E \xi(h-h^*)(X)\right| \leq \rho \|h-h^*\|_{L_2}^2;
$$
thus
$$
P_N {\cal L}^{{\cal K}}_h \geq \E {\cal L}^{{\cal K}}_h - 3\rho\|h-h^*\|_{L_2}^2.
$$
otherwise, if $\|h-h^*\|_{L_2} \leq r_{{U}}$,
\begin{align*}
P_N {\cal L}^{{\cal K}}_h \geq & \E {\cal L}^{{\cal K}}_h - 2\left|\frac{1}{N} \sum_{i=1}^N \xi_i (h-h^*)(X_i) - \E \xi(h-h^*)(X)\right| - \|h-h^*\|_{L_2}^2
\\
\geq & \E{\cal L}_h^{{\cal K}} -3 r_{{U}}^2,
\end{align*}
as claimed.
\endproof

\begin{Lemma} \label{lemma:structure-hat-F}
For a sample ${\cal D}=(X_i,Y_i)_{i=1}^N$, let $d={\rm diam}(V({\cal D}),L_2)$ and recall that $f^*={\rm argmin}_{f \in F} \E(f(X)-Y)^2$. If ${\cal D} \in {\cal A}_{f^*}$ then
\begin{description}
\item{1.} $f^* \in V({\cal D})$, and
\item{2.} for every $v \in V({\cal D})$, $\|v-Y\|_{L_2}^2 \leq \|f^*-Y\|_{L_2}^2+6\max\left\{r_{{U}}^2, d^2/400\right\}$.
\end{description}
\end{Lemma}
\proof
Fix ${\cal D} =(X_i,Y_i)_{i=1}^N \in {\cal A}_{f^*}$. Recall that if $f_1,f_2 \in F \subset {U}$ and $\|f_1-f_2\|_{L_2} \geq r_{{U}}$, one has
$$
\alpha \|f_1-f_2\|_{L_2} \leq {\rm Med}_\ell(|f_1-f_2|(X_i))_{i=1}^N \leq \beta \|f_1-f_2\|_{L_2}.
$$
In addition, applying Lemma \ref{lemma:key-observation} for ${\cal K}=F$, it follows that for every $f \in F$,
\begin{equation} \label{eq:non-standard-oracle-in-proof}
0 \leq \E {\cal L}^F_f \leq P_N {\cal L}^F_f  + 3\max\{r_{{U}}^2, \rho \|f-f^*\|_{L_2}^2\}.
\end{equation}
 Let $\hat{f}$ be the empirical minimizer in $F$ and consider the following two cases: if $\|\hat{f}-f^*\|_{L_2} \geq r_{{U}}$ then ${\rm Med}_\ell(|\hat{f}-f^*|(X_i))_{i=1}^N \geq \alpha \|\hat{f}-f^*\|_{L_2}$; alternatively, $\|\hat{f}-f^*\|_{L_2} \leq r_{{U}}$. Therefore, by \eqref{eq:non-standard-oracle-in-proof}
\begin{align*}
P_N {\cal L}^F_{\hat{f}} = & P_N (\hat{f}-Y)^2 -P_N(f^*-Y)^2 \geq -3\max\{r_{{U}}^2, \rho\|\hat{f}-f^*\|_{L_2}^2\}
\\
\geq & -3\max\left\{r_{{U}}^2, \rho\alpha^{-2} {\rm Med}_\ell^2\left(|\hat{f}-f^*|(X_i)\right)_{i=1}^N\right\},
\end{align*}
implying that $f^* \in V({\cal D})$.

Turing to the second part, note that $\hat{f} \in V({\cal D})$, $P_N {\cal L}_{\hat f}^F \leq 0$ and that for every $u \in U$,
$$
{\rm Med}_\ell(|\hat{f}-u|(X_i))_{i=1}^N \leq \beta \max\{r_{{U}}, \|\hat{f}-u\|_{L_2}\} \leq \beta \max\{r_{{U}}, d\}.
$$
Hence, it follows from the definition of $V({\cal D})$ that for every $v \in V({\cal D})$,
\begin{align} \label{eq:in-proof-5.4}
P_N {\cal L}^{F}_v \leq &  P_N {\cal L}^{F}_{\hat{f}} + 3\max\left\{ r_{{U}}^2, \rho\alpha^{-2}{\rm Med}_\ell^2\left(|\hat{f}-v|(X_i)\right)_{i=1}^N\right\} \nonumber
\\
\leq & 3  \max\left\{ r_{{U}}^2, \rho \left(\frac{\beta}{\alpha}\right)^2d^2\right\}.
\end{align}
Combining \eqref{eq:in-proof-5.4} with \eqref{eq:non-standard-oracle-in-proof}, for every $v \in V({\cal D})$,
\begin{equation*}
\|v-Y\|_{L_2}^2 - \|f^*-Y\|_{L_2}^2 =  \E {\cal L}^{F}_v \leq
6\max\left\{r_{{U}}^2, d^2/400\right\},
\end{equation*}
by the choice of $\rho$.
\endproof

\noindent {\bf Proof of Theorem \ref{thm:procedure}.}
One has to show that the assumptions of Lemma \ref{lemma:V-class} hold for $V({\cal D}_1)$ and for $W({\cal D}_1)$ for the sample ${\cal D}_2$. By Lemma \ref{lemma:structure-hat-F}, $f^* \in V({\cal D}_1)$ and thus, for every $v \in V({\cal D}_1)$,
$$
\|v-Y\|_{L_2}^2 - \|f^*-Y\|_{L_2}^2 \leq \max\left\{6r_{{U}}^2, d^2/50\right\}.
$$
Also, applying Lemma \ref{lemma:key-observation} for $W({\cal D}_1) \equiv W \subset {U}$, it follows that if ${\cal D}_2 \in A_{w^*}$ then for every $w \in W$,
\begin{align*}
\E {\cal L}_w^W \leq & P_N {\cal L}_w^W + 3\max\{r_{U}^2, \rho \|w-w^*\|^2_{L_2}\}
\\
\leq & P_N {\cal L}_w^W + \max\{3r_{U}^2, \|w-w^*\|^2_{L_2}/50\}
\end{align*}
where $P_N$ is the empirical mean relative to ${\cal D}_2$. Thus, the assumptions of Lemma \ref{lemma:V-class} are verified, completing the proof of Theorem \ref{thm:procedure}.
\endproof

\section{The events ${\cal A}_{u_0}$} \label{sec:the-event}
The final part of the of the proof of Theorem \ref{thm:main} focuses on the events ${\cal A}_{u_0}$. We will show that for every $u_0 \in {U}$, ${\cal A}_{u_0}$ is a high probability event provided that $\alpha,\beta$ and $\ell$ are properly chosen constants that depend only on $q$ and $L$, and that $r_U^2=r_{{\rm opt}}$ for the right choice of constants.

\subsection{An almost isometric lower estimate} \label{sec:isometric}
The main result of this section is an `almost isometric' lower bound on $\inf_{f \in {\cal F}} \frac{1}{N}\sum_{i=1}^N f^2(X_i)$ for an arbitrary class ${\cal F}$.

The small-ball method, introduced in \cite{Men-LWC,Men-GAFA,MenKolt,Men-LWCG}, may be used to show that if ${\cal F}$ satisfies the small-ball condition with constants $\kappa_0$ and $\eps$ then
$$
\inf_{\{f \in {\cal F}, \|f\|_{L_2} \geq r\}} \frac{1}{N}\sum_{i=1}^N \left(\frac{f(X_i)}{\|f\|_{L_2}}\right)^2 \geq c_0,
$$
but $c_0=c_0(\kappa_0,\eps)$ is a constant that {\it need not be close to $1$}. To obtain an almost isometric result rather than an `isomorphic' one, a slightly stronger assumption is required.

\begin{Theorem} \label{thm:isometric}
For every $2<q \leq 4$ and $L \geq 1$ there exist constants $c_1$ and $c_2$ that depend only on $q$ and $L$ for which the following holds. Let ${\cal K}={\rm star}({\cal F})$ and assume that for every $h \in {\cal K}-{\cal K}$, $\|h\|_{L_q} \leq L \|h\|_{L_2}$. Set $\gamma_1=q/2(q-1)$ and $\gamma_2=(q-2)/2(q-1)$, and let $0<\zeta<1$ and $r$ for which
\begin{equation} \label{eq:fixed-points-in-isometric}
\E\|G\|_{({\cal K}-{\cal K}) \cap rD} \leq  \zeta \sqrt{N}r, \ \
\ \
\E\sup_{h \in ({\cal K}-{\cal K}) \cap rD} \frac{1}{\sqrt{N}}\left|\sum_{i=1}^N \eps_i h(X_i) \right| \leq \zeta \sqrt{N}r.
\end{equation}
Then, with probability at least $1-2N\exp(-c_1 \zeta^{\gamma_1}N)$, if $h \in {\cal K}-{\cal K}$ and $\|h\|_{L_2} \geq r$,
\begin{equation} \label{eq:isometric-lower-in-thm}
\frac{1}{N}\sum_{i=1}^N h^2(X_i) \geq \|h\|_{L_2}^2 \left(1-c_2 \zeta^{\gamma_2}\right).
\end{equation}
\end{Theorem}

It is highly likely that the exponents $\gamma_1$ and $\gamma_2$ are not optimal. For example, when $q=4$ one would expect an estimate of $(1-c_2\zeta^{1/2})\|h\|^2_{L_2}$ rather than $(1-c_2\zeta^{1/3})\|h\|^2_{L_2}$ that follows from Theorem \ref{thm:isometric}. Fortunately, this gap has little effect on the proof of Theorem \ref{thm:main}: once the value of $\alpha$ and $\beta$ is chosen, $\rho=(\alpha/20\beta)^2$ and $\zeta$ satisfies $\rho=c_2\zeta^{\gamma_2}$; thus $\zeta$ is a small but fixed constant that depends only on $q$ and $L$, and the suboptimal power in \eqref{eq:isometric-lower-in-thm} will be of little significance in what follows.
\vskip0.4cm
For the proof of Theorem \ref{thm:isometric} we will first present an almost isometric lower estimate for a finite set. In the general case, that set will be an appropriate net which approximates ${\cal F}$, and the final step in the proof will be an upper estimate on the empirical `approximation errors'.

\vskip0.5cm
\subsubsection{An estimate for a single function}
Given integers $N$ and $m$ and a function $f \in L_q$ for some $2<q \leq 4$, set
\begin{equation*}
\phi(f) =
\begin{cases}
f & \ \ {\rm if} \ \ |f| \leq \left(\frac{N}{m}\right)^{1/q} \|f\|_{L_q}
\\
\\
\left(\frac{N}{m}\right)^{1/q}{\rm sgn}(f) & \ \ {\rm if} \ \ |f| > \left(\frac{N}{m}\right)^{1/q} \|f\|_{L_q}.
\end{cases}
\end{equation*}
Hence, $\phi$ is a truncation of $f$ at a level that is selected according to the $L_q$ space to which $f$ belongs, the sample size $N$ and a parameter $m$ that will be used to calibrate the probability estimate.

Observe that pointwise, $|\phi(f)| \leq |f|$, and given a sample $(X_i)_{i=1}^N$, set
$$
I_f = \left\{i: \left(\phi(f)\right)(X_i) = f(X_i)\right\}=\left\{i: |f(X_i)| \leq (N/m)^{1/q}\|f\|_{L_q}\right\}.
$$

\begin{Theorem} \label{thm:decomposition-single-function}
There exist absolute constants $c_0$, $c_1,c_2$ and for $2<q \leq 4$ there are constants $c_3,c_4$ that depend only on $q$ for which the following holds. If $N \geq c_0m$, then with probability at least $1-2\exp(-c_1m)$,
\begin{description}
\item{1.} $|I_f| \geq N -c_2m$, and
\item{2.} for every $J \subset \{1,...,N\}$ with $|J| \leq 4 m$,
\begin{align*}
& \left(\|f\|_{L_2}^2-c_3\|f\|_{L_q}^2\left(\frac{m}{N}\right)^{1-(2/q)}\right)  \leq  \frac{1}{N} \sum_{i \in I_f \backslash J} f^2(X_i)
\\
\leq & \left(\|f\|_{L_2}^2+c_4\|f\|_{L_q}^2\left(\frac{m}{N}\right)^{1-(2/q)}\right).
\end{align*}
\end{description}
\end{Theorem}

\proof
Observe that $Pr(|f| \geq (N/m)^{1/q}\|f\|_{L_q}) \leq {m}/{N}$. Using a standard binomial estimate applied to the event $\{|f| \geq (N/m)\|f\|_{L_q}\}$, it follows that for $0<u<N/4m$, $$
Pr(|I_f^c| \geq um) \leq \binom{N}{um} \cdot \left(\frac{m}{N}\right)^{um} \leq \left(\frac{e}{u}\right)^{um},
$$
and the first claim follows.

Next, consider $h=(\phi(f))^2$. Since $\|h\|_{L_\infty} \leq (N/m)^{2/q} \|f\|_{L_q}^2$ and $|\phi(f)| \leq |f|$,
\begin{align*}
\|h\|_{L_2}^2 = & \E (\phi(f))^q \cdot (\phi(f))^{4-q} \leq \|f\|_{L_q}^q \cdot \left(\frac{N}{m}\right)^{-1+4/q} \|f\|_{L_q}^{4-q}
\\
= & \left(\frac{N}{m}\right)^{-1+4/q} \cdot \|f\|_{L_q}^4.
\end{align*}
Hence, by Bernstein's inequality (see, e.g., \cite{VW}) for $h=(\phi(f))^2$,
$$
\left|\frac{1}{N} \sum_{i=1}^N \left(\phi(f)\right)^2(X_i) - \E \left(\phi(f)\right)^2 \right| \leq \left(\frac{m}{N}\right)^{1-2/q}\|f\|_{L_q}^2
$$
with probability at least
$$
1-2\exp\left(-c_2N \min\left\{\frac{u^2}{\|h\|_{L_2}^2},\frac{u}{\|h\|_{L_\infty}}\right\}\right)=1-2\exp(-c_2m).
$$
Recall that if $i \in I_f^c$ then $|(\phi(f))(X_i)| = (N/m)^{1/q} \|f\|_{L_q}$, and by the first part of the claim, $|I_f^c| \leq c_3m$. Therefore,
\begin{equation*}
\sum_{i \in I_f^c} \left(\phi(f)\right)^2(X_i) \leq \frac{|I_f^c|}{N} \cdot \left(\frac{N}{m}\right)^{2/q} \|f\|_{L_q}^2
\leq c_4\left(\frac{m}{N}\right)^{1-2/q} \|f\|_{L_q}^2.
\end{equation*}
Also,
$$
\E \left(\phi(f)\right)^2 \geq \E f^2 \IND_{\{|f| \leq (N/m)^{1/q}\|f\|_{L_q}\}} = \E f^2 - \E f^2 \IND_{\{|f| >  (N/m)^{1/q}\|f\|_{L_q}\}},
$$
and
\begin{align*}
& \E f^2 \IND_{\{|f| >  (N/m)^{1/q}\|f\|_{L_q}\}} = \int_0^\infty 2t Pr \left(|f|\IND_{\{|f| > (N/m)^{1/q}\|f\|_{L_q}\}} > t \right)dt
\\
\leq & \left(\frac{N}{m}\right)^{2/q}\|f\|_{L_q}^2 Pr (|f| > (N/m)^{1/q}\|f\|_{L_q}) + \int_{(N/m)^{1/q}\|f\|_{L_q}}^\infty 2t Pr(|f|>t) dt
\\
\leq & \frac{q}{q-2} \left(\frac{m}{N}\right)^{1-2/q}\|f\|_{L_q}^2;
\end{align*}
thus,
$$
\E (\phi(f))^2 \geq \E f^2 - \frac{q}{q-2}\left(\frac{m}{N}\right)^{1-2/q}\|f\|_{L_q}^2.
$$
Combining these observations, with probability at least $1-2\exp(-c_2m)$,
\begin{align*}
\frac{1}{N}\sum_{i \in I_f} f^2(X_i) = &\frac{1}{N}\sum_{i \in I_f} (\phi(f))^2(X_i) =
\frac{1}{N}\sum_{i=1}^N (\phi(f))^2(X_i)-\frac{1}{N}\sum_{i \in I_f^c} (\phi(f))^2(X_i)
\\
\geq & \E (\phi(f))^2 - c_4 \left(\frac{m}{N}\right)^{1-2/q}\|f\|_{L_q}^2
\\
\geq & \E f^2 -c(q) \left(\frac{m}{N}\right)^{1-2/q}\|f\|_{L_q}^2,
\end{align*}
and, in a similar fashion,
$$
\frac{1}{N}\sum_{i \in I_f} f^2(X_i) \leq \E f^2 + c(q) \left(\frac{m}{N}\right)^{1-2/q}\|f\|_{L_q}^2
$$
for a constant $c(q)$ that depends only on $q$.

Finally, note that if $i \in I_f$ then $|f(X_i)| \leq (N/m)^{1/q}\|f\|_{L_q}$. Hence, for every $J \subset I_f$,
$$
\frac{1}{N}\sum_{j \in J} f^2(X_i) \leq \frac{|J|}{N} \cdot \left(\frac{N}{m}\right)^{2/q}\|f\|_{L_q}^2,
$$
which completes the proof.
\endproof

\subsubsection{A uniform lower bound}
Let ${\cal F} \subset L_2$ and $\eta>0$, and set
$$
\Phi_N({\cal F},\eta)=\E \sup_{h \in ({\cal F} - {\cal F}) \cap \eta D} \left|\frac{1}{N}\sum_{i=1}^N \eps_i h(X_i) \right|+\eta.
$$
When the underlying class ${\cal F}$ or the sample size $N$ are clear, we will abuse notation and write $\Phi(\eta)$ instead of $\Phi_N({\cal F},\eta)$.
\vskip0.4cm
Let $N(\eps,{\cal F},L_2)$ be the minimal number of open $\eps$-balls with respect to the $L_2$ norm that are needed to cover ${\cal F}$, and set
$$
e_m({\cal F})=\inf\{ \eps>0 : \log N(\eps,{\cal F},L_2) \leq m\}
$$
to be the $m$-th entropy number of ${\cal F}$. The centres of the balls are called a {\it minimal cover of ${\cal F}$}.

\begin{Lemma} \label{lemma:lower-global-preliminary}
Let ${\cal F} \subset rS(L_2)$, set $2<q \leq 4$ and assume that for every $f_1,f_2 \in {\cal F} \cup \{0\}$, $\|f_1-f_2\|_{L_q} \leq L \|f_1-f_2\|_{L_2}$. If $1 \leq m \leq N$ and $e_m({\cal F}) \leq \eta $ then with probability at least $1-2N\exp(-c_1m)$,
$$
\inf_{f \in {\cal F}} \frac{1}{N}\sum_{i=1}^N f^2(X_i) \geq r^2\left(1-c_2(q,L)\left(\left(\frac{N}{m}\right)^{1/2}\frac{\Phi(\eta)}{r}+
\left(\frac{m}{N}\right)^{1-2/q}\right)\right),
$$
where $c_1$ is an absolute constant and $c_2$ is a constant that depend only on $L$ and $q$.
\end{Lemma}

\proof
Fix an integer $m$ and let ${\cal F}^\prime$ be a minimal $\eta$-cover of ${\cal F}$ with respect to the $L_2$ norm. Since $\eta \geq e_m({\cal F})$ it follows that $\log |{\cal F}^\prime| \leq m$.

Let $f \in {\cal F}$ and set $\pi f \in {\cal F}^\prime$ for which $\|f -\pi f\|_{L_2} \leq \eta$. Put $v_j = \pi f(X_j)$ and $u_j=(f-\pi f)(X_j)$, and observe that if $I \subset \{1,...,N\}$ then
\begin{align*}
\sum_{i=1}^N f^2(X_i) \geq & \sum_{i \in I} f^2(X_i) =  \sum_{i \in I} \left(\pi f (X_i) + (f-\pi f)(X_i)\right)^2
\\
\geq & \sum_{i \in I} v_i^2 - 2 \left(\sum_{i \in I} v_i^2 \right)^{1/2} \cdot \left(\sum_{i \in I} u_i^2 \right)^{1/2}.
\end{align*}

Let $I_{\pi f}=\{i : \pi f(X_i)=(\phi(\pi f))(X_i)\}$ be as in Theorem \ref{thm:decomposition-single-function}; set $J_f \subset \{1,...,N\}$ to be the union of the set of the largest $2 m$ coordinates of $\left(|f-\pi f|(X_i)\right)_{i=1}^N=(|u_i|)_{i=1}^N$ and the set of the largest $2m$ coordinates of $(|\pi f(X_i)|)_{i=1}^N=(v_i)_{i=1}^N$. Applying Theorem \ref{thm:decomposition-single-function}, the union bound and the $L_q$-$L_2$ norm equivalence, there is an absolute constant $c_1$ and a constant $c_2$ that depends only on $q$ for which, with probability at least $1-2\exp(-c_1m)$, for every $\pi f \in {\cal F}^\prime$,
$$
\left(1-c_2 L\left(\frac{m}{N}\right)^{1-(2/q)}\right) \leq \frac{1}{N} \sum_{i \in I_{\pi f} \backslash J_f } \left(\frac{\pi f(X_i)}{\|\pi f\|_{L_2}}\right)^2 \leq  \left(1+c_2L\left(\frac{m}{N}\right)^{1-(2/q)}\right).
$$

Next, one has to obtain a high probability estimate on the `coordinate distribution' of the vector $\left(|u_i|\right)_{i=1}^N$. To that end, fix $t>0$ and observe that by symmetrization and contraction arguments (see, e.g. \cite{LT,VW}),
\begin{align*}
& t \E \sup_{f \in {\cal F}} |\{i : |f-\pi f|(X_i) \geq t\}| \leq \E \sup_{f \in {\cal F}} \sum_{i=1}^N |f-\pi f|(X_i)
\\
\leq & 2\E \sup_{f \in {\cal F}} \left|\sum_{i=1}^N \eps_i (f-\pi f)(X_i)\right| + N\sup_{f \in {\cal F}} \E |f-\pi f | \leq 2 N \Phi(\eta).
\end{align*}
Fix $t_j$ to be named later and apply Talagrand's concentration inequality for bounded empirical processes \cite{Tal94,Led,BLM} to the class of indicator functions $\{\IND_{\{|f-\pi f| \geq t_j\}} : f \in {\cal F}\}$.
Thus, with probability at least $1-2\exp(-m)$, for every $f \in {\cal F}$,
\begin{equation*}
\frac{1}{N} \sum_{i=1}^N \IND_{\{|f-\pi f | \geq t_j\}}(X_i) \leq c_3 \left( \E\sup_{f \in {\cal F}} \frac{1}{N} \sum_{i=1}^N \IND_{\{|f-\pi f | \geq t_j\}}(X_i) + \sqrt{\frac{m}{N}} \sigma_j + \frac{m}{N}\right)=(*)_j,
\end{equation*}
where $\sigma_j = \sup_{f \in {\cal F}} Pr^{1/2}( |f-\pi f| \geq t_j)$. By the $L_q$ and $L_2$ norm equivalence,
$$
\sigma_j^2 \leq \sup_{f \in {\cal F}} \frac{\E|f-\pi f|^q}{t_j^q} \leq \left(\frac{L \eta}{t_j}\right)^q.
$$
Therefore, if $j \geq 2m$ and $t_j=c_4(q,L)\Phi(\eta)N/j$,
$$
(*)_j \leq c_3\left( \frac{2\Phi(\eta)}{t_j} + \sqrt{\frac{m}{N}} \cdot \left(\frac{L\eta}{t_j}\right)^{q/2} + \frac{m}{N}\right) \leq \frac{j}{N},
$$
because $\Phi(\eta) \geq \eta$.

Summing for $2m \leq j \leq N$, with probability at least $1-2N\exp(-m)$, for every $j \geq 2m$ and every $f \in {\cal F}$,
$$
\left|\left\{i : |(f-\pi f)(X_i)| \geq c_4\Phi(\eta) \frac{N}{j}\right\}\right| \leq j.
$$
And, on that event for $j \geq 2m$,
$$
u_j^* \leq c_4\Phi(\eta) \frac{N}{j},
$$
where $(u_i^*)_{i=1}^N$ denotes a non-increasing rearrangement of $(|u_i|)_{i=1}^N$.

Hence,
$$
\left(\sum_{j >2m} (u_j^*)^2 \right)^{1/2} \leq c_4 \Phi(\eta) N \left(\sum_{j >2m} j^{-2}\right)^{1/2} \leq c_5(q,L) \Phi(\eta)N/\sqrt{m},
$$
and setting $I=I_{\pi f} \backslash J_f$,
$$
\left(\sum_{i \in I} v_i^2\right)^{1/2} \cdot \left(\sum_{i \in I} u_i^2\right)^{1/2} \leq c_6(q,L) N^{3/2} r \Phi(\eta)/\sqrt{m}.
$$
Recalling the lower estimate on $\sum_{i \in I} v_i^2$, it follows that
$$
\frac{1}{N}\sum_{i \in I} f^2(X_i) \geq r^2 \cdot \left(1-c(L,q)\left(\left(\frac{m}{N}\right)^{1-(2/q)}+
\left(\frac{\Phi(\eta)}{r}\right)\cdot\left(\frac{N}{m}\right)^{1/2}\right)\right).
$$
\endproof

\noindent{\bf Proof of Theorem \ref{thm:isometric}.} Let $0<\zeta<1$ and set $r$ for which \eqref{eq:fixed-points-in-isometric} holds. Recall that ${\cal K}={\rm star}({\cal F})$, and thus ${\cal K}-{\cal K}$ is star-shaped around $0$. Hence, it is standard to verify that if $r^\prime \geq r$ than
$$
\E \sup_{w \in ({\cal K}-{\cal K}) \cap r^\prime D} \left|\frac{1}{N}\sum_{i=1}^N \eps_i w(X_i) \right| \leq \zeta r^\prime
$$
and
$$
\E\|G\|_{({\cal K}-{\cal K}) \cap r^\prime D} \leq \zeta \sqrt{N}r^\prime.
$$
Consider the class ${\cal F}_r={\rm star}({\cal F}) \cap r S(L_2)={\cal K} \cap rS(L_2)$. Set $\eta =c \E\|G\|_{{\cal F}_r}/\sqrt{m}$ and note that by Sudakov's minoration (see, e.g. \cite{Pis,LT}), $\eta \geq e_m({\cal F}_r)$, provided that $c$ is a well-chosen absolute constant. Therefore,
\begin{equation*}
\Phi_N({\cal F}_r,\eta) \leq \E \sup_{w \in ({\cal K}-{\cal K}) \cap 2r D} \left|\frac{1}{N}\sum_{i=1}^N \eps_i w(X_i) \right| + \eta \leq 2\zeta r + \frac{c\E \|G\|_{{\cal F}_r}}{\sqrt{m}}.
\end{equation*}
Set $m =\theta N$ for a constant $0<\theta<1$ to be specified later; thus,
$$
\frac{\Phi_N({\cal F}_r,\eta)}{r} \leq 2\zeta + \frac{c}{\sqrt{\theta}} \cdot \frac{\E\|G\|_{{\cal F}_r}}{r \sqrt{N}}.
$$
Moreover, ${\cal F}_r \subset ({\cal K}-{\cal K}) \cap rD$. By the choice of $r$, $\E\|G\|_{{\cal F}_r} \leq \E\|G\|_{({\cal K}-{\cal K}) \cap rD} \leq \zeta r \sqrt{N}$, and
$$
\frac{\Phi_N({\cal F}_r,\eta)}{r} \leq \frac{c_1\zeta}{\sqrt{\theta}}.
$$
Thanks to Lemma \ref{lemma:lower-global-preliminary}, with probability at least $1-2N\exp(-c_2 \theta N)$, for every $f \in {\cal F}_r$,
\begin{align} \label{eq:lower-in-proof}
\nonumber \frac{1}{N}\sum_{i=1}^N f^2(X_i) \geq & r^2 \left(1-c(q,L)\left(\left(\frac{m}{N}\right)^{1-2/q}+
\left(\frac{N}{m}\right)^{1/2}\frac{\Phi(\eta)}{r}\right)\right)
\\
\geq & r^2\left(1-c(q,L)\left(\theta^{1-2/q}+\frac{c\zeta}{\theta} \right) \right).
\end{align}
Setting $\theta=\zeta^{q/2(q-1)}$, the claim follows for $f \in {\cal F}_r={\rm star}({\cal F}) \cap r S(L_2)$. 

Finally, since \eqref{eq:lower-in-proof} is positive homogeneous and ${\rm star}({\cal F})$ is star-shaped around $0$, it also holds on the same event when $f \in {\rm star}({\cal F})$ and $\|f\|_{L_2} >r$.
\endproof

\subsection{The Median of means as a crude measure of distances} \label{sec:median}
As noted above, the results of \cite{Men-LWC,Men-GAFA,MenKolt,Men-LWCG} show that the small-ball method suffices to ensure that with probability at least $1-2\exp(-cN)$,
$$
\alpha^2 \|f-h\|_{L_2}^2 \leq \frac{1}{N} \sum_{i=1}^N (f-h)^2(X_i)
$$
for well chosen constants $\alpha$ and $c$ that depend only on the small-ball condition in ${\cal F}$, and for every $f,h \in {\cal F}$ whose $L_2$-distance is not `too small'. However, if class members do not have well-behaved tails, the probability that
$$
\alpha^2 \|f-h\|_{L_2}^2 \leq \frac{1}{N} \sum_{i=1}^N (f-h)^2(X_i) \leq \beta^2 \|f-h\|_{L_2}^2
$$
even for a {\it single} pair $f,h \in {\cal F}$ may be rather small; certainly not of the order of $1-2\exp(-cN)$.  Unfortunately, this means that the empirical mean is a poor two-sided estimator of distances, as it lacks stability: if $f-h$ is a heavy-tailed function, there will be at least one very large value of $|(f-h)(X_i)|$, and that will destroy any hope of having
$$
\frac{1}{N} \sum_{i=1}^N (f-h)^2(X_i) \leq \beta^2 \|f-h\|_{L_2}^2
$$
unless $\beta$ is very large.

To bypass this obstacle, we will use the more stable median-of-means functional.
\vskip0.3cm
Let us begin by showing that very little `mixing' is needed for an empirical mean to satisfy a small-ball estimate with a rather high (constant) probability.
\begin{Lemma} \label{lemma:Berry-Esseen}
For every $q>2$ and $L \geq 1$ there are constants $\ell$ and $\kappa_0$ that depend only on $q$ and $L$ for which the following holds. If $\|Z\|_{L_q} \leq L \|Z\|_{L_2}$ and $Z_1,...,Z_\ell$ are independent copies of $Z$, then
$$
Pr\left(\frac{1}{\ell} \sum_{i=1}^\ell |Z_i| \geq \kappa_0 \|Z\|_{L_2} \right) \geq \frac{3}{4}.
$$
\end{Lemma}

\proof Since $\|Z\|_{L_q} \leq L \|Z\|_{L_2}$, it follows from a standard application of the Paley-Zygmund inequality (see, e.g., \cite{dlPG}) that $Z$ satisfies a small-ball condition with constants $c_1$ and $c_2$ that depend only of $q$ and $L$. Therefore,
\begin{equation} \label{eq:L1-L2-mpment-in-proof}
\|Z\|_{L_1} \geq c_1\|Z\|_{L_2} Pr(|Z| \geq c_1\|Z\|_{L_2}) \geq c_1c_2\|Z\|_{L_2}.
\end{equation}
By an appropriate version of the Berry-Esseen inequality for independent copies of $Z \in L_q$ for $q>2$ \cite{MOO}, if $\ell \geq c_3(q,L)$ then
$$
\sup_{t \in \R} \left| Pr \left( \frac{1}{\ell} \sum_{i=1}^\ell |Z_i| \geq  \E|Z|+\frac{t\||Z|-\E|Z|\|_{L_2}}{\sqrt{\ell}} \right) - Pr\left(g \geq t\right) \right| \leq 0.05.
$$
Take $t < 0$ to be the largest for which $Pr\left(g \geq t\right) \geq 0.8$. Applying \eqref{eq:L1-L2-mpment-in-proof}, if $\ell \geq 4t^2/(c_1c_2)^2$ then
$\E |Z| \geq 2|t|\|Z\|_{L_2}/\sqrt{\ell}$, and
$$
\E|Z|+\frac{t\||Z|-\E|Z|\|_{L_2}}{\sqrt{\ell}} \geq (c_1c_2/2)\|Z\|_{L_2}.
$$
 Therefore, setting $\kappa_0=c_1c_2/2$ (which depends only on $q$ and $L$),
$$
Pr\left(\frac{1}{\ell}\sum_{i=1}^N |Z_i| \geq \kappa_0\|Z\|_{L_2}\right) \geq \frac{3}{4}.
$$
\endproof

Fix $2<q \leq 4$ and $L \geq 1$, and set $\ell$ and $\kappa_0$ as in Lemma \ref{lemma:Berry-Esseen}.
Without loss of generality, assume that $N=\ell M$ for an integer $M$, and recall that for $v \in \R^N$, ${\rm Med}_\ell(v)$ is the median of the vector of means performed in the $M$ blocks $I_0,...,I_{M-1}$.

\begin{Theorem} \label{thm:median}
For every $2<q \leq 4$ and $L \geq 1$ there exists constants $c_1,c_2,c_3$ and $\alpha<1<\beta$ that depend only on $q$ and $L$, for which the following holds. Let ${\cal F} \subset L_2$, put ${\cal K}={\rm star}({\cal F})$ and assume that for every $w \in {\cal K}$, $\|w\|_{L_q} \leq L \|w\|_{L_2}$. Set $r>0$ that satisfies
\begin{equation*}
E\|G\|_{({\cal K}-{\cal K}) \cap rD} \leq c_1 \sqrt{N}r, \ \ \
\E\sup_{h \in ({\cal K}-{\cal K}) \cap rD} \frac{1}{\sqrt{N}}\left|\sum_{i=1}^N \eps_i h(X_i) \right| \leq c_2 \sqrt{N} r.
\end{equation*}
Then, with probability at least $1-2\exp(-c_3 N)$, for every $w \in {\cal K}$ for which $\|w\|_{L_2} \geq r$,
$$
\alpha \|w\|_{L_2} \leq {\rm Med}_\ell\left(|w(X_i)|\right)_{i=1}^N \leq \beta\|w\|_{L_2}.
$$
Moreover, on the same event, if $\|w\|_{L_2} \leq r$ then
$$
{\rm Med}_\ell\left(|w(X_i)|\right)_{i=1}^N \leq \beta r.
$$
\end{Theorem}

The proof of Theorem \ref{thm:median} follows the same lines as the proof of Theorem 4.3 from \cite{Men-LWCG}. It is based on the following observation.

\begin{Lemma} \label{lemma:single-variable}
There are absolute constants $c_1$ and $c_2$ for which the following holds. Consider $Z \in L_2$ that satisfies a small-ball condition with constants $\kappa_0$ and $\eps$. If $Z_1,...,Z_N$ are independent copies of $Z$, then with probability at least $1-2\exp(-c_1\delta^2\eps N)$ there is a subset $I \subset \{1,...,N\}$, $|I| \geq (1-\delta)\eps N$, and for every $i \in I$,
$$
\kappa_0 \|Z\|_{L_2} \leq |Z_i| \leq c_2\|Z\|_{L_2}/\sqrt{\delta \eps}.
$$
\end{Lemma}

\proof Fix $0<\delta<1$ and let $A=\{\kappa_0 \|Z\|_{L_2} \leq |Z| \leq 3\|Z\|_{L_2}/\sqrt{\delta\eps}\}$. Combining the small-ball condition and Chebyshev's inequality, $Pr(A) \geq 1-(1+\delta/3)\eps$. Let $\eta$ to be a selector (i.e., a $\{0,1\}$-valued random variable) with mean $(1+\delta/3)\eps$ and set $\eta_1,...,\eta_N$ to be independent copies of $\eta$. A standard concentration argument shows that with probability at least $1-2\exp(-c_1\delta^2\eps N)$,
\begin{equation*}
|\{i: \eta_i=1\}|= \sum_{i=1}^N \eta_i \leq (1+\delta/3)^2 \eps N \leq (1+\delta)\eps N,
\end{equation*}
and the claim follows.
\endproof

\noindent{\bf Proof of Theorem \ref{thm:median}.} Let $\ell$ and $\kappa_0$ be as in Lemma \ref{lemma:Berry-Esseen} and recall that the two constants depend only on $q$ and $L$. Assume, without loss of generality, that $M=N/\ell$ is an integer, set $\eps=3/4$ and fix $0<\delta<1$ for which $(1-\delta)\eps=0.6$. Set ${\cal K}={\rm star}({\cal F})$, let $\zeta_1$ and $\zeta_2$ to be named later and put $r>0$ that satisfies
\begin{align*}
& E\|G\|_{{\cal K} \cap rD} \leq \zeta_1 \sqrt{N}r \ \ \ {\rm and}
\\
& \E\sup_{h \in ({\cal K}-{\cal K}) \cap rD} \frac{1}{\sqrt{N}}\left|\sum_{i=1}^N \eps_i h(X_i) \right| \leq \zeta_2 \sqrt{N}r.
\end{align*}
Let $v \in {\cal K}$ and set
$$
{\cal M}_v=\frac{1}{\ell} \sum_{i=1}^\ell |v|(X_i);
$$
thus $({\cal M}_{v,j})_{j=0}^{M-1} = \left(\ell^{-1}\sum_{i \in I_j} |v|(X_i)\right)_{j=0}^{M-1}$
are $M$ independent copies of the random variable ${\cal M}_v$, which, by Lemma \ref{lemma:Berry-Esseen} satisfies the small-ball condition with constants $\kappa_0$ and $\eps=3/4$.

One may verify that
$$
\frac{3}{4}\kappa_0 \|v\|_{L_2} \leq \|{\cal M}_v\|_{L_2} \leq \|v\|_{L_2}.
$$
By Lemma \ref{lemma:single-variable}, with probability at least $1-2\exp(-c_1\delta^2 \eps M)=1-2\exp(-c_2N)$,
there is $J \subset \{0,...,M-1\}$, $|J| \geq  (1-\delta/2)\eps M$, and for every $j \in J$,
$$
\frac{3}{4}\kappa_0^2 \|v\|_{L_2} \leq {\cal M}_{v,j} \leq c_3  \|v\|_{L_2}/\sqrt{\delta \eps} =c_4\|v\|_{L_2}.
$$
Hence, the same assertion holds uniformly for $\exp(c_2 N/2)$ random variables of the form ${\cal M}_v$. And, in particular, for every $v \in {\cal V}_r \subset {\cal K} \cap r S(L_2) \equiv {\cal K}_r$, which is a maximal $\eta$-separated set for a choice of $\eta$ large enough to ensure that $|{\cal V}_r| \leq \exp(c_2 N/2)$.

Therefore, with probability at least $1-2\exp(-c_2 N/2)$, for every $v \in {\cal V}_r$ there is a subset $J_v \subset \{0,...,M-1\}$ of cardinality at least $(1-\delta/2)\eps M$ and for every $j \in J_v$,
\begin{equation} \label{eq:mean-medians-est}
\frac{3}{4}\kappa_0^2 r = \frac{3}{4}\kappa_0^2 \|v\|_{L_2} \leq {\cal M}_{v,j} \leq c_4 r.
\end{equation}
By Sudakov's inequality applied to the set ${\cal K}_r$ and using the choice of $r$, one may select
$$
\eta = c_5 \frac{\E \|G\|_{{\cal K}_r}}{\sqrt{c_2 N /2}} \leq c_6 \zeta_1r.
$$
Next, consider the empirical oscillation term: for every $f \in {\cal K}_r$, let $\pi f $ be the best approximation with respect to the $L_2$ distance of $f$ in ${\cal V}_r$. Set $u_f=\IND_{\{|f-\pi f| > 3\kappa_0^2r/8\}}$, consider the class of indicator functions ${\cal U}_r=\{ u_f : f \in {\cal K}_r \}$ and let
$$
\psi(X_1,...,X_N)=\sup_{u_f \in {\cal U}_r} \frac{1}{N} \sum_{i=1}^N u_f(X_i).
$$
 By the bounded differences inequality (see, for example, \cite{BLM}), with probability at least $1-\exp(-c_7t^2)$,
$$
\psi(X_1,...,X_N) \leq \E \psi + \frac{t}{\sqrt{N}}.
$$
To estimate $\E \psi$ from above, set $\phi(t)=t/(3\kappa_0^2 r/8)$. Observe that for every $u_f \in {\cal U}_r$, $u_f(X) \leq \phi(|f-\pi f|(X))$, and that by the Gin\'{e}-Zinn symmetrization theorem \cite{GZ84,LT} and the choice of $r$,
\begin{align*}
& \E \sup_{u_f \in {\cal U}_r} \frac{1}{N} \sum_{i=1}^N u_f(X_i) \leq \E \sup_{f \in {\cal K}_r} \frac{1}{N}\sum_{i=1}^N \phi(|f-\pi f|(X_i))
\\
\leq & \E \sup_{f \in {\cal K}_r} \left|\frac{1}{N}\sum_{i=1}^N \phi(|f-\pi f|(X_i)) - \E \phi(|f-\pi f|(X_i)) \right| + \sup_{f \in {\cal K}_r} \E \phi(|f-\pi f|)
\\
\lesssim & \frac{1}{\kappa_0^2 r} \cdot \left( \E \sup_{f \in {\cal K}_r} \left|\frac{1}{N} \sum_{i=1}^N \eps_i (f-\pi f)(X_i) \right|+\sup_{f \in {\cal K}_r} \|f-\pi f\|_{L_2} \right)
\\
\lesssim & \frac{1}{\kappa_0^2 r} \cdot \left(\zeta_2 r + \eta \right) \leq \frac{\delta \eps}{4 \ell}
\end{align*}
when $\zeta_1 \sim \kappa_0^2/\ell $ and $\zeta_2 \sim \kappa_0^2/\ell $, and thus depend only on $q$ and $L$.

Setting $t=\delta\eps\sqrt{N}/4\ell$, it follows that with probability at least $1-2\exp(-c_8(q,L)N)$,  for every $f \in {\cal K}_r={\rm star}({\cal F}) \cap r S(L_2)$,
$$
|\{i : |f-\pi f|(X_i) \geq (3\kappa_0^2/8)r\}| \leq \delta \eps N/2\ell=\delta \eps M/2.
$$
Therefore, at most $\delta \eps M/2$ of the $M$ `bins' $I_j$ contain a sample point $X_i$ for which $|f-\pi f|(X_i) \geq (3\kappa_0^2/8)r$; on the remaining $(1-\delta \eps/2)M$ bins,
$$
\frac{1}{\ell} \sum_{i \in I_j} |f-\pi f|(X_i) \leq (3\kappa_0^2/8)r.
$$
Hence, with probability at least $1-2\exp(-c_9(q,L)N)$, for every $f \in {\rm star}({\cal F}) \cap rS(L_2)$ there is a subset of $\{0,...,M-1\}$ of cardinality at least $(1-\delta)\eps M = 0.6M$, on which
\begin{align} \label{eq:signs}
\nonumber
|{\cal M}_{f,j}| \geq & |{\cal M}_{\pi f,j}| - |{\cal M}_{f-\pi f,j}| \geq (3\kappa_0^2/4) r - (3\kappa_0^2/8)r =(3\kappa_0^2/8)\|f\|_{L_2}, \ \ \ {\rm and}
\\
|{\cal M}_{f,j}| \leq & |{\cal M}_{\pi f,j}| + |{\cal M}_{f-\pi f,j}| \leq (c_4+ 3\kappa_0^2/8)\|f\|_{L_2}. \end{align}
Moreover, since the estimates are positive homogeneous and ${\rm star}({\cal F})$ is star-shaped around $0$, \eqref{eq:signs} is true on the same event when $f \in {\rm star}({\cal F})$ and $\|f\|_{L_2} \geq r$.
The claim follows by recalling that $\ell$ and $\kappa_0$ depend only on $q$ and $L$, and selecting $0<\alpha<3\kappa_0/8$ and $\beta \geq c_4 + 3\kappa_0^2/8$.

The proof of the second part is almost identical: ${\cal V}_r$ is defined exactly as above, and for every $f \in {\rm star}({\cal F}) \cap r D$, $\pi f$ is the best approximation in ${\cal V}_r$; thus, $\|f-\pi f\|_{L_2} \leq 2r$. Just as in the proof of the first part, with probability at least $1-2\exp(-c_9(q,L)N)$,  for every $f \in {\rm star}({\cal F}) \cap r D$,
$$
|\{i : |f-\pi f|(X_i) \geq 3r\}| \leq \delta \eps M/2.
$$
Thus, on at least $(1-\delta)\eps M = 0.6M$ of the `bins'
$$
|{\cal M}_{\pi f,j}|+|{\cal M}_{f-\pi f,j}| \leq c_{10}(q,L)r,
$$
and one may choose $\beta=\max\{c_4+3\kappa_0^2/8,c_{10}\}$.
\endproof

\section{Proof of Theorem \ref{thm:main}}

Observe that the second the third conditions in the definition of ${\cal A}_{u_0}$ are independent of $u_0$, and we shall begin by verifying those.
\vskip0.4cm
Given the base class $F$, recall that $U=\{(f_1+f_2)/2 : f_1,f_2 \in F\}$ and that $H={\rm star}(U-U)$. Thus, for every $h \in H$, $\|h\|_{L_q} \leq L \|h\|_{L_2}$. Let $r_0$ be the infimum of the set of all $r>0$ for which
\begin{align} \label{eq:in-proof-1-9}
& \E\|G\|_{H \cap rD} \leq  \zeta_1(q,L)\sqrt{N}r \ \ \ {\rm and} \nonumber
\\
& \E \sup_{h \in (H-H) \cap rD}\left|\frac{1}{\sqrt{N}}\sum_{i=1}^N \eps_i h(X_i)\right| \leq \zeta_2(q,L)\sqrt{N}r
\end{align}
for constants $\zeta_1$ and $\zeta_2$ as in Theorem \ref{thm:median}.

Since $H$ and $H-H$ are star-shaped around $0$, \eqref{eq:in-proof-1-9} holds for every $r \geq r_0$. Invoking Theorem \ref{thm:median} for ${\cal F}=U-U$ and $r=2r_0$, there are constants $\alpha \leq 1 \leq \beta$ and $\ell$ that depend only on $q$ and $L$ for which, with probability at least $1-2\exp(-c_0(q,L)N)$, for every $h \in H$,
\begin{description}
\item{$\bullet$} if $\|h\|_{L_2} \geq 2r_0$ then $\alpha \|h\|_{L_2} \leq {\rm Med}_\ell(|h(X_i)|)_{i=1}^N \leq \beta\|h\|_{L_2}$;
\item{$\bullet$} if $\|h\|_{L_2} \leq 2r_0$ then ${\rm Med}_\ell(|h(X_i)|)_{i=1}^N \leq \beta \cdot 2r_0$.
\end{description}
In particular, for any $r_U \geq 2r_0$, the third condition in the definition of ${\cal A}_{u_0}$ is verified.

Next, let $\alpha$ and $\beta$ be as above and set $\rho=(\alpha/20\beta)^2$. Consider Theorem \ref{thm:isometric} for ${\cal F}=U-U$ (and in which case, ${\cal K}={\rm star}(U-U)=H$). Recall that $\gamma_1=q/2(q-1)$ and $\gamma_2=(q-2)/2(q-1)$ and set $\zeta_3$ by $\rho \sim \zeta_3^{\gamma_2}$, and in particular, $\zeta_3$ depends only on $q$ and $L$. Set $r_1$ for which
\begin{align*}
& \E\|G\|_{(H-H) \cap r_1D} \leq \zeta_3 \sqrt{N}r_1, \ \ {\rm and}
\\
& \E\sup_{h \in (H-H)\cap r_1D} \left|\frac{1}{\sqrt{N}}\sum_{i=1}^N \eps_i h(X_i)\right| \leq \zeta_3 \sqrt{N}r_1.
\end{align*}
It follows that with probability at least
$$
1-2N\exp(-c_1\zeta_3^{\gamma_1}N) \geq 1-2\exp(-c_2(q,L)\rho^{\gamma_1/\gamma_2}N),
$$
if $h \in H-H$ and satisfies $\|h\|_{L_2} \geq r_1$ then
$$
\frac{1}{N}\sum_{i=1}^N h^2(X_i) \geq (1-\rho)\|h\|_{L_2}^2.
$$
Moreover, since $0 \in H$, the same is true for every difference $h=u_1-u_2 \in U-U \subset H-H$ provided that $\|u_1-u_2\|_{L_2} \geq r_1$. Thus, the second part in the definition of ${\cal A}_{u_0}$ holds for $r_{U} \geq r_1$.
\vskip0.4cm
Turning to the first part of the definition of ${\cal A}_{u_0}$ (which does depend on $u_0$), one may
apply Lemma \ref{lemma:multi-and-zero-level-sym} to the set $U$ and for
$$
r_2 \geq  r_M(F,\rho/4,\delta/2,u_0).
$$
Setting $\xi=u_0(X)-Y$ and $\xi_i = u_0(X_i)-Y_i$, it follows that with probability at least $1-\delta$, for every $u \in U$,
$$
\left|\frac{1}{N}\sum_{i=1}^N \xi_i (u-u_0)(X_i) - \E \xi (u-u_0)\right| \leq \rho \max\left\{\|u-u_0\|_{L_2}^2,r_2^2\right\},
$$
as required.

Finally, one may combine all the above conditions, by noting that for $c_3=\rho/4$, $c_4=\min\{\zeta_2,\zeta_3\}$ and $c_5=\min\{\zeta_1,\zeta_3\}$, the choice of $r^2=r_{{\rm opt}}(F,\delta,c_3,c_4,c_5)$ is a valid choice in all of the above. Hence, for every $u_0 \in U$,
$$
Pr({\cal A}_{u_0}) \geq 1-\delta-2\exp(-c_6(q,L)N),
$$
and Theorem \ref{thm:main} follows from Theorem \ref{thm:procedure}.
\endproof

\subsection{Proof of Corollary \ref{cor:subgaussian}} \label{sec:finite-dictionary}
Let $F$ be a finite dictionary. While a  learning procedure can only guarantee an error rate of the order of $\sqrt{N^{-1}\log M}$, one may show that the aggregation procedure suggested above leads to a much better estimate.

Let us begin by reformulating Corollary \ref{cor:subgaussian}:
\begin{Theorem} \label{thm:finite-dictionary}
For every $L \geq 1$ and $q>2$ there exist a constant $c_1$ that depends only on $L$ and $q$ for which the following holds. Let $F=\{f_1,...,f_M\}$ and assume that for $w \in {\rm span}(F)$ and every $p \geq 2$, $\|w\|_{L_p} \leq L\sqrt{p}\|w\|_{L_2}$. Assume further that $Y \in L_q$ for some $q>2$. Then for every $0<\delta<1$, with probability at least $1-\delta$,
$$
\E\left((\tilde{f}(X)-Y)^2|(X_i,Y_i)_{i=1}^N\right) \leq \E(f^*(X)-Y)^2 + c_1 \delta^{-2/q}\log(2/\delta) \|f^*-Y\|^2_{L_q}\frac{\log M}{N}.
$$
\end{Theorem}

As all the assumptions of Theorem \ref{thm:main} are satisfied here, what is left is to identify $r_{{\rm opt}}$. To that end, note that $|U-U| \leq M^4$. Thus, for every $r>0$ there is $k_r \leq M^8$ and functions $(w_{i,r})_{i=1}^{k_r}$, that satisfy $\|w_{i,r}\|_{L_2} \leq r$ and
$$
{\rm star}(H-H) \cap rD \subset \left\{ \lambda w_{i,r} : 1 \leq i \leq k_r, 0 \leq \lambda \leq 1\right\} \equiv W_r.
$$
By the moment equivalence in ${\rm span}(F)$, a straightforward chaining argument and the Majorizing Measures Theorem (see, e.g., \cite{MPT} for similar arguments) it follows that
$$
\E \sup_{w \in W_r} \left|\frac{1}{\sqrt{N}} \sum_{i=1}^N \eps_i w(X_i)\right| \leq c_1L \E\|G\|_{W_r}.
$$
And, it is standard to verify that
$$
\E\|G\|_{W_r} = \E \sup_{1 \leq i \leq k_r} G_{w_i} \leq c_2 r \sqrt{\log M}.
$$
Therefore, if $N \geq c_3(L,\zeta) \log M$, then
$$
r_{Q,1}(F,\zeta)=r_{Q,2}(F,\zeta)=0.
$$
Turning our attention to $r_M$, one may invoke the following fact from \cite{Men-Bern}:
\begin{Theorem} \label{thm:multiplier}
Let $\xi \in L_q$ for some $q>2$ and assume that for every $f,h \in {\cal F} \cup \{0\}$ and every $p \geq 2$, $\|f-h\|_{L_p} \leq L \sqrt{p}\|f-h\|_{L_2}$. Then, for every $u,w >1$, with probability at least
\begin{equation} \label{eq:prob-est-in-multiplier}
1-c_0(q)w^{-q}\frac{\log^q N}{N^{q/2-1}} - 2\exp\left(-c_1(L)u^2 \left(\frac{\E\|G\|_{{\cal F}}}{{\rm diam}({\cal F},L_2)}\right)^2\right),
\end{equation}
one has
$$
\sup_{f \in {\cal F}} \left|\frac{1}{\sqrt{N}}\sum_{i=1}^N \eps_i \xi_i f(X_i)\right| \leq c_2(q)Lwu\|\xi\|_{L_q} \E\|G\|_{\cal F}.
$$
\end{Theorem}
\vskip0.4cm
For every $u_0 \in U$ let ${\cal F}={\rm star}(U-u_0) \cap rD$. Since $|U-u_0|=|U| \leq M^2$ then by Theorem \ref{thm:multiplier}, and with probability as in \eqref{eq:prob-est-in-multiplier},
$$
\sup_{f \in {\cal F}} \left|\frac{1}{\sqrt{N}}\sum_{i=1}^N \eps_i \xi_i f(X_i)\right| \leq c(q)Lwu\|\xi\|_{L_q}r \sqrt{\log M}.
$$
Therefore, if
$$
r \geq c(q)\frac{Lwu}{\zeta} \cdot \|\xi\|_{L_q}\sqrt{\frac{\log M}{N}},
$$
then
$$
\sup_{f \in {\cal F}} \left|\frac{1}{\sqrt{N}}\sum_{i=1}^N \eps_i \xi_i f(X_i)\right| \leq \zeta \sqrt{N} r^2.
$$
Clearly, for any nontrivial class ${\cal F}$, $\E\|G\|_{\cal F} \gtrsim {\rm diam}({\cal F},L_2)$; thus, setting $w \sim (1/\delta)^{1/q}$ and $u \sim \sqrt{\log(2/\delta)}$, with probability at least $1-\delta$
$$
r_M \leq c_1(q) \frac{L}{\zeta} \cdot \delta^{-1/q}\log^{1/2}(2/\delta) \|\xi\|_{L_q}\sqrt{\frac{\log M}{N}},
$$
which completed the proof of Theorem \ref{thm:finite-dictionary}.

\subsection{A remark on the bounded case}
Let us briefly mention a way in which one may obtain a version of Theorem \ref{thm:main} when both the dictionary and the target are assumed to be bounded in $L_\infty$, but $F$ may be infinite.

As noted in \cite{LM-lin-agg}, an $L_\infty$ type of assumption is of a very different nature than an assumption on norm equivalence: the former does not lead to a useful small-ball estimate on class members, and in particular, the proofs presented in Section \ref{sec:the-event} do not hold in that case.

Fortunately, there are highly potent tools at one's disposal when bounded classes are concerned, namely, Talagrand's concentration inequality for bounded empirical processes and the contraction principle for empirical and Bernoulli processes indexed by bounded classes (see, e.g., \cite{LT,VW,BLM}). Using that well established machinery, one may show that ${\cal A}_{u_0}$ is a high probability event. In fact, thanks to the two-sided concentration estimates, the argument is much simpler.

For example, assuming that the functions involved are bounded by $1$ almost surely and applying a contraction argument, it follows that with high probability and in expectation,
\begin{align*}
& \sup_{u \in U} \left|\frac{1}{N}\sum_{i=1}^N (u_0(X_i)-Y_i)(u-u_0)(X_i) - \E(u_0(X)-Y)(u-u_0)(X) \right|
\\
\lesssim & \sup_{u \in U} \left|\frac{1}{N}\sum_{i=1}^N \eps_i (u-u_0)(X_i)\right|
\end{align*}
and
\begin{equation*}
\sup_{u \in U} \left|\frac{1}{N}\sum_{i=1}^N (u-u_0)^2(X_i)-\E(u-u_0)^2\right|
\lesssim \sup_{u \in U} \left|\frac{1}{N}\sum_{i=1}^N \eps_i (u-u_0)(X_i)\right|,
\end{equation*}
where $(\eps_i)_{i=1}^N$ are independent, symmetric $\{-1,1\}$-valued random variables that are independent of $(X_i,Y_i)_{i=1}^N$.

Moreover, the multiplier and quadratic processes concentrate well around their mean, leading to a natural complexity parameter that is rather similar to $r_{{\rm opt}}$, and to an exponential probability estimate.

The obvious downside in this concentration-contraction based argument is that it totally eliminates the dependence on the distance between $F$ and $Y$ (see the discussion in \cite{Men-LWC,Men-LWCG} for more details). As an outcome, the estimate in the bounded case does not improve when the problem becomes more `realizable'.

\footnotesize {

\end{document}